\documentclass[11pt]{amsart}
\usepackage{amsfonts,amssymb}
\usepackage[english]{babel}
\usepackage{amscd}
\usepackage{hyperref}
\hypersetup{colorlinks,citecolor=blue}

\newcommand{\BZ}{{\mathbf{Z}}}

\newcommand{\BR}{{\mathbf{R}}}
\newcommand{\BC}{{\mathbf{C}}}
\newcommand{\BQ}{{\mathbf{Q}}}
\newcommand{\BP}{{\mathbf{P}}}
\newcommand{\cH}{{\mathcal{H}}}

\newcommand{\gd}{\delta}

\newcommand{\gs}{\sigma}

\newcommand{\gep}{\varepsilon}

\newcommand{\ro}{\varrho}

\newcommand{\gf}{\varphi}

\newcommand{\mC}{{\mathcal{C}}}
\newcommand{\mD}{{\mathcal{D}}}

\def\e{\varepsilon}
\def\hy{-\allowhyphens}
\newcommand{\ti}[1]{\tilde{#1}}

\newcommand{\ol}[1]{\overline{#1}}
\newcommand{\sol}[1]{{\,\overline{\!#1}\!}}

\newcommand{\diam}{{\rm diam}}
\newcommand{\diag}{{\rm diag}}
\newcommand{\SL}{\mathbf{SL}}
\newcommand{\GL}{\mathbf{GL}}

\newcommand{\PGL}{\mathbf{PGL}}

\newcommand{\Sp}{\mathbf{Sp}}

\newcommand{\Fix}{\text{Fix}}
\newcommand{\Isom}{{\rm Isom}}
\newcommand{\Aff}{{\rm Aff}}
\newcommand{\Orth}{\mathbf{O}}

\gdef\overto#1{{\buildrel{#1}\over\longrightarrow}}

\newcommand{\ball}{\mathrm{B}}
\newcommand{\sphere}{\mathrm{S}}

\def\dist{\text{dist}}
\def\Aut{\text{Aut}}

\def\se{\subseteq}

\newcommand{\ip}[2]{\langle{#1},{#2}\rangle}

\newcommand{\Set}[2]{\left\{ {#1}\,:\, {#2}\right\}  }

\def\no{n\raise4pt\hbox{\tiny o}\kern+.2em}

\numberwithin{equation}{section}

\newtheorem{mthm}{Theorem}

\newtheorem{prop}{Proposition}[section]
\newtheorem{thm}[prop]{Theorem}
\newtheorem{lem}[prop]{Lemma}
\newtheorem{cor}[prop]{Corollary}
\newtheorem{conj}[prop]{Conjecture}
\theoremstyle{definition}

\newtheorem{defn}[prop]{Definition}
\newtheorem{rem}[prop]{Remark}
\newtheorem{rems}[prop]{Remarks}
\newtheorem{exam}[prop]{Example}

\newtheorem{ques}[prop]{Question}
\newtheorem{ob}[prop]{Observation}
\newtheorem{clm}[prop]{Claim}

\newtheorem{facts}[prop]{Facts}

\begin{document}
\author[U. Bader]{Uri~Bader}
\address{U.B. Technion -- Israel Institute of Technology, Haifa, Israel}
\author[A. Furman]{Alex~Furman}
\address{A.F. University of Illinois at Chicago, Chicago, USA}
\author[T. Gelander]{Tsachik~Gelander}
\address{T.G. Yale university, USA}
\curraddr{The Hebrew university, Jerusalem, Israel}
\author[N. Monod]{Nicolas~Monod}
\address{N.M. University of Chicago, Chicago, USA}
\curraddr{Universit\'e de Gen\`eve, Switzerland}
\date{\today}
\thanks{U.B. partially supported by ISF grant~100146, 
A.F. partially supported by NSF grants DMS-0094245
and DMS 0604611, T.G. partially supported by NSF grant
DMS-0404557 and BSF grant 2004010, 
N.M. partially supported by FNS (CH) and NSF (US)}

\title[Property~$(T)$ for actions on Banach spaces]{Property~$(T)$ and rigidity\\ for actions on Banach spaces}

\begin{abstract}
We study property~$(T)$ and the fixed point property for actions
on $L^p$ and other Banach spaces. We show that property~$(T)$
holds when $L^2$ is replaced by $L^p$ (and even a
subspace/quotient of $L^p$), and that in fact it is independent of
$1\leq p<\infty$. We show that the fixed point property for $L^p$
follows from property~$(T)$ when $1<p< 2+\e$. For simple Lie
groups and their lattices, we prove that the fixed point property
for $L^p$ holds for any $1< p<\infty$ if and only if the rank is
at least two. Finally, we obtain a superrigidity result for
actions of irreducible lattices in products of general groups on
superreflexive spaces.
\end{abstract}

\maketitle
\let\languagename\relax 
\section{Introduction and the Main Results}

\subsection{}
Since its introduction by Kazhdan in~\cite{Kazhdan},
property~$(T)$ became a fundamental concept in mathematics with a
wide range of applications to such areas as:

\smallskip

\noindent \textbullet\ The structure of infinite groups~-- finite
generation and finite abelanization of higher rank
lattices~\cite{Kazhdan}, obstruction to free or amalgamated
splittings~\cite{Wa}, \cite{Al}, \cite{Ma-amlg} structure of
normal subgroups~\cite{Ma-normal} etc.;

\noindent
\textbullet\ Combinatorics~-- the first construction of expanders~\cite{Ma-expanders} (see~\cite{Lub});

\noindent \textbullet\ Operator algebras~-- factors of type II$_1$
whose fundamental group is countable~\cite{ConnesT} or even
trivial~\cite{Popa04}); rigidity theorems for the factors
associated to Kazhdan group~\cite{Popa1};

\noindent \textbullet\ Ergodic theory~-- rigidity results related
to Orbit Equivalence~\cite{Popa2},\cite{Hjorth}; the
Banach--Ruziewicz problem~\cite{Ma-mean},\cite{Sul};

\noindent \textbullet\ Smooth dynamics~-- local
rigidity~\cite{FM1},\cite{FM2}; actions on the
circle~\cite{Navas1}
(and~\cite{Pressley-Segal},\cite{Reznikov_AT}).

\smallskip

It has also been an important tool in providing interesting
(counter) examples: to Day's ``von Neumann
conjecture''~\cite[5.6]{Gro-hyp} and in the context of the
Baum--Connes conjecture~\cite{HLS} (related to \cite{Gro-random}).

\medskip

Initially defined in terms of unitary representations,
property~$(T)$ turned out to be equivalent to Serre's
property~$(FH)$ -- a fixed point property for affine isometric
actions on Hilbert spaces that can be rephrazed as cohomological
vanishing. (The equivalence holds for $\sigma$-compact groups, in
particular all locally compact second countable groups, and was
proved by Delorme~\cite{Delorme} and Guichardet~\cite{Guichardet}.
As pointed out by Y.~de Cornulier, uncountable discrete groups
that have G.~Bergman's cofinality property~\cite{Bergman04}
have~$(FH)$ but fail~$(T)$.) Some of the above applications use
this latter characterization. Recently Shalom~\cite{Shalom}
described the \emph{reduced} 1-cohomology with unitary
coefficients for irreducible lattices in products of completely
general locally compact groups. This led to a list of new rigidity
results and added such lattices to the list of ``naturally rigid''
groups. For further details and more references on these topics,
we suggest the monography~\cite{HV} and the
forthcoming~\cite{BHV}.

\subsection{}
Motivated by these broad themes:
property~$(T)$, property~$(FH)$, lattices in semisimple groups and in
general products, we study similar notions in the broader framework of Banach spaces rather than
Hilbert spaces.  Some of the results below apply to general
\emph{superreflexive} Banach spaces, whilst some are specific to
the subclass of $L^{p}(\mu)$\hy spaces with $1<p<\infty$.
(A Banach space is superreflexive if it admits an equivalent uniformly convex norm,
see Proposition~\ref{P:superreflex-ucus} below.)

One of the motivations to consider such questions came from the work of
Fisher and Margulis~\cite{FM1}, \cite{FM2}, in which an $L^{p}$ analogue
of property~$(T)$ with $p\gg2$ allowed them to weaken smoothness
assumptions in their results.

The harder question of fixed point results for \emph{affine actions} on $L^p$ for $p\gg2$ (see Theorem~B below)
has applications e.g. for actions on the circle~\cite{Navas2}, \cite{BHV}.

\subsection{}
Let $G$ be a topological group and $B$ a Banach space. By a
\emph{linear isometric $G$\hy representation} on $B$, we shall
mean a continuous homomorphism $\ro:G\to\Orth(B)$ where $\Orth(B)$
denotes the (``orthogonal'') group of all invertible linear
isometries $B\to B$ (see Lemma~\ref{L:cont-reps} for a
clarification of the continuity assumption).
We say that such a representation \emph{almost has invariant vectors} if
\begin{equation}\label{e:almostinv}
    \forall\, \text{compact subset } K\se G,\qquad \inf_{\|v\|=1} \diam(\ro(K)v)=0.
\end{equation}
Denote by $B^{\ro(G)}$ the closed subspace of $G$\hy fixed vectors;
the $G$\hy representation $\ro$ descends to a linear
isometric $G$\hy representation $\ro^{\prime}$ on $B^{\prime}=B/B^{\ro(G)}$
(see Remark~\ref{R:TB} for more details in the case of superreflexive spaces).
We shall use the following as a Banach space analogue of Kazhdan's property~$(T)$:
\begin{defn}\label{D:TB}
Let $B$ be a Banach space. A topological group $G$ is said to have property~$(T_B)$
if for any continuous linear isometric $G$\hy representation $\ro:G\to\Orth(B)$ the quotient
$G$\hy representation
$\ro^{\prime}:G\to \Orth(B/B^{\ro(G)})$ does not almost have $G$\hy invariant vectors.
\end{defn}
Note that if $B$ is a Hilbert space, $\ro^{\prime}$ is isomorphic
to the restriction of $\ro$ to the orthogonal complement
$(B^{\ro(G)})^{\perp}$ of the subspace of $\ro(G)$\hy invariants.
Thus for Hilbert spaces the above definition agrees with Kazhdan's
property~$(T)$.

Let $\mu$ be a $\sigma$\hy finite measure on a standard Borel space $(X,\mathcal{B})$.
We are most interested in the family $L^{p}(\mu)$, $1<p<\infty$, of Banach spaces,
which are close relatives of Hilbert spaces.
They also possess a rich group of linear isometries $\Orth(L^{p}(\mu))$.

\begin{mthm} \label{T:TLp}
Let $G$ be a locally compact second countable group.
If $G$ has Kazhdan's property~$(T)$ then $G$ has property~$(T_B)$
for Banach spaces $B$ of the following types:
\begin{itemize}
\item[{\rm (i)}]
        $L^{p}(\mu)$ for any $\sigma$\hy finite measure $\mu$ and any $1\leq p<\infty$.
\item[{\rm (ii)}]
        A closed subspace of $L^{p}(\mu)$ for any $1<p<\infty$, $p\neq 4,6,8,\ldots$.
\item[{\rm (iii)}]
        A quotient space of $L^{p}(\mu)$ for any $1<p<\infty$, $p\neq
                \frac{4}{3},\frac{6}{5},\frac{8}{7},\ldots$.
\end{itemize}
If $G$ has~$(T_{L^{p}([0,1])})$ for some $1<p<\infty$ then
$G$ has Kazhdan's property~$(T)$.
\end{mthm}

\subsection{}
Next we consider group actions by isometries on Banach spaces. By
the Mazur--Ulam theorem, such actions are always affine with the
linear part being isometric as well (working with \emph{real}
Banach spaces).

\begin{defn} \label{D:FB}
We say that $G$ has property~$(F_B)$ if any continuous action of $G$ on $B$ by affine isometries
has a $G$\hy fixed point.
\end{defn}

When $B$ is a Hilbert space this is precisely Serre's
property~$(FH)$. Delorme~\cite{Delorme} and
Guichardet~\cite{Guichardet} proved that properties (T) and~$(FH)$
are equivalent for $\sigma$-compact groups.
Below we summarize the relations between properties $(T)$
and $(F_{B})$ which hold for general groups.
\begin{thm}\label{T:T-FB-general}\addtocounter{footnote}{1}
For a locally compact second countable group $G$ we have
\begin{enumerate}
\item
    $(F_{B})$ implies $(T_{B})$ for any Banach space $B$.
\item
    $(T)$ implies $(F_{B})$ for closed subspaces $B$ of $L^{p}(\mu)$ where $1<p\le 2$.\\
Likewise for subspaces of $L^1$ and of the pseudo-normed spaces $L^p(\mu)$,\\
$0<p< 1$, except one obtains only bounded orbits instead of fixe
points\footnote{See Example~\ref{E:L10} for an example without
fixed point.}.

\item
    $(T)$ also implies $(F_{B})$ for closed subspaces of $L^{p}(\mu)$
    for $2\le p<2+\e$, where $\e=\e(G)>0$ might depend
    on the Kazhdan group $G$.
\end{enumerate}
\end{thm}
\begin{rems}
(1) is essentially due to Guichardet~\cite{Guichardet} as his proof of
$(FH)\Rightarrow(T)$ applies to all Banach spaces.
We give two proofs for (2) reducing the problem, in both, to
one of the proofs of $(T)\Rightarrow(FH)$.
We note that the particular case of $p=1$ in (2) is one of the results of~\cite{RS}.
Statement (3) is due to Fisher and Margulis (unpublished).
With their kind permission we have included their argument
here (see Section~\ref{S:FM}).
\end{rems}
\medskip

The above results imply that any locally compact group $G$
with Kazhdan's property~$(T)$ has property~$(T_{L^{p}})$ for all $1<p<\infty$,
and has the fixed point property~$(F_{L^{p}})$ for $1<p<2+\e(G)$.
It turns out, however, that many Kazhdan groups (e.g. hyperbolic ones)
do not have property~$(F_{L^{p}})$ for large values of $p$.

Indeed, in his study of $L^{p}$\hy cohomology, Pansu~\cite{Pansu}
proved that $\Sp_{n,1}(\BR)$ and cocompact lattices in these groups have a
non-trivial first $L^p$\hy cohomology $L^{p}H^{1}$ for all $p>4n+2$.
This is equivalent to asserting that for $p>4n+2$ these groups admit
fixed-point-free affine isometric actions on $L^{p}(G)$ with linear part being the
regular representation.
Hence these groups do not have property~$(F_{L^{p}})$ for $p>4n+2$,
whilst enjoying $(T)$.

More generally, $L^{p}H^{1}(\Gamma)$ and hence $H^{1}(\Gamma,\ell^{p}\Gamma)$
is non-zero for any non-elementary hyperbolic group
when $p$ is large enough. Indeed,
Bourdon and Pajot identify this
cohomology with a Besov space of functions on the boundary, which
they prove to be non-trivial as soon as $p$ is larger than the
Hausdorff dimension of an Ahlfors-regular metric on the boundary,
see Corollaire~6.2 in~\cite{Bourdon-Pajot}.
Again, this contradicts $(F_{L^{p}})$ for large $p$.

More recently, using Mineyev's homological
bicombings~\cite{Mineyev}, Yu~\cite{Yu} gave a very short proof
that any hyperbolic group $\Gamma$ admits a \emph{proper} action
by affine isometries on $\ell^{p}(\Gamma\times\Gamma)$ if $p$ is
large enough. This is a strong negation of $(F_{L^{p}})$ for
hyperbolic groups and all their infinite subgroups. The
corresponding strenghtening of the above mentioned~\cite{Pansu},
\cite{Bourdon-Pajot} for rank one Lie (or algebraic) groups $G$
has been established by Cornulier--Tessera--Valette
in~\cite{Cornulier-Tessera-Valette}: For any $p>1$ larger than the
Hausdorff dimension of the boundary, there is a \emph{proper}
affine isometric action on $L^p(G)$ whose linear part is the
regular representation. In particular, this holds for
$\Sp_{n,1}(\BR)$ when $p>4n+2$.

\subsection{}
Our next goal is now, by contrast, to establish $(F_{L^{p}})$ for certain groups.
It is often remarked that property~$(T)$ for (simple) higher rank Lie groups
and their lattices is more robust than property~$(T)$ enjoyed by
the rank one groups $\Sp_{n,1}(\BR)$ and many other Gromov hyperbolic groups.
In view of the preceding discussion of hyperbolic groups and $\Sp_{n,1}(\BR)$, the following result might
be viewed as yet another evidence supporting this view.

\begin{mthm}\label{T:higherrankFLp}
Let $G=\prod_{i=1}^{m} \mathbf{G}_{i}(k_{i})$, where $k_{i}$ are
local fields (of any characteristic), $\mathbf{G}_{i}(k_{i})$ are
$k_{i}$\hy points of Zariski connected simple $k_{i}$\hy algebraic
groups $\mathbf{G}_{i}$. Assume that each simple factor
$\mathbf{G}_{i}(k_{i})$ has $k_{i}$\hy rank $\ge 2$.

\nobreak
Then $G$ and the lattices in $G$ have property~$(F_{B})$ for all
$L^{p}(\mu)$\hy related spaces $B$ as in {\rm (i)--(iii)} in
Theorem~\ref{T:TLp}, assuming $1<p<\infty$.
\end{mthm}

\subsection{}
A broader class of spaces in which we propose to study properties
$(T_{B})$ and $(F_{B})$ consists of \emph{superreflexive} spaces,
which can be defined as topological vector spaces isomorphic to
uniformly convex Banach spaces\footnote{For spaces that are only
\emph{strictly convex}, the fixed point property always
fails~\cite{Brown-Guentner},\cite{Haagerup-Przybyszewska}.}. In
this context we consider linear representations (resp. affine
actions) which are \emph{uniformly equicontinuous}; more
concretely, for any given norm compatible with the topology, the
class of all such linear representations (resp. affine actions) is
that of \emph{uniformly bounded} linear representations (resp.
\emph{uniformly Lipschitz} affine actions). It turns out that such
representations (resp. actions) can always be viewed as isometric
with respect to some equivalent norm that is \emph{simultaneously}
uniformly convex and uniformly smooth
(Proposition~\ref{P:superreflex-affine}).

Note that whether a given linear $G$-representation almost contains invariant vectors
or not, in the sense of~\eqref{e:almostinv}, does not depend on a particular norm among
all mutually equivalent norms. Hence we can make the following
\begin{defn}\label{D:superTFB}
Let $B$ be a superreflexive topological vector space and $G$ a locally compact second countable group.
We say that $G$ has property~$(\sol{T}_{B})$ if for every uniformly equicontinuous linear
representation $\ro$ of $G$ on $B$ the quotient $G$-representation on $B/B^{\ro(G)}$
does not almost have invariant vectors.

Likewise, $G$ has $(\sol{F}_{B})$ if every uniformly equicontinuous affine
$G$-action on $B$ has a fixed point.
\end{defn}
\begin{conj} \label{C:higherrankFB}
Higher rank groups $G=\prod \mathbf{G}_{i}(k_{i})$
as in Theorem~\ref{T:higherrankFLp} and their lattices have property~$(\sol{F}_{B})$,
and hence $(\sol{T}_{B})$, for all superreflexive $B$.
\end{conj}
\begin{rem}
To support this conjecture let us point out the following:

(1) Much of our proof of Theorem~\ref{T:higherrankFLp} is done in the broad
context of uniformly equicontinuous affine actions on general superreflexive spaces
except for one argument~-- a version of relative property~$(T_{B})$, whose proof is special to $L^{p}$\hy related
spaces.

(2) V. Lafforgue proved~\cite{La} that the group $\PGL_3(\BQ_p)$
has property~$(\sol{T}_B)$ for all superreflexive $B$ (his result
is actually stronger, in that he allows linear representations
with \emph{slowly growing}, rather \emph{uniformly bounded}
Lipschitz norms, see Theorem 3.2, Definition 0.2 and the
discussion preceding it in \cite{La}). Combined with our proof of
Theorem~\ref{T:higherrankFLp} it implies for example that
$\SL_{n}(\BQ_{p})$, $n\ge 4$, has property~$(\sol{F}_{B})$.

(3) Y. Shalom has proved (unpublished) that for Hilbert spaces $\mathcal H$
higher rank groups (and their lattices) have property~$(\sol{F}_{\mathcal H})$, and hence
$(\sol{T}_{\mathcal H})$, whilst rank one groups
have neither $(\sol{F}_{\mathcal H})$ nor $(\sol{T}_{\mathcal H})$.
\end{rem}

\subsection{}
One way to generalize the context of semisimple (non-simple)
Lie/algebraic groups is simply to consider general products
$G=G_1\times\cdots\times G_n$ of $n\geq2$ \emph{arbitrary}
topological groups. In the absence of any assumption on the
factors $G_i$, one can still establish splitting results for
uniformly equicontinuous affine $G$\hy actions on superreflexive
spaces.
\begin{mthm}\label{thm_splitting}%
Let $G=G_1\times \cdots\times G_n$ be a product of topological groups with
a continuous action by uniformly equicontinuous affine maps on a superreflexive topological
vector space $B$ without $G$\hy fixed point.
Assume that the associated linear $G$\hy representation $\ro$ does not almost have
non-zero invariant vectors.

\nobreak Then there is a $G$\hy invariant closed complemented
affine subspace $\underline{B}\se B$ and an affine equicontinuous
$G$\hy equivariant isomorphism $\underline{B}\cong B_1 \oplus
\cdots \oplus B_n$, where each $B_i$ is a superreflexive Banach
space with an equicontinuous affine $G$\hy action factoring
through $G\to G_i$.
\end{mthm}
\begin{rems}\label{R:thm_splitting}

(1)~If $G$ has property~$(\sol{T}_B)$ then the assumption that $\ro$
does not almost have invariant vectors is redundant.

(2)~In the particular case where $B$ is a Hilbert space and $G$
locally compact acting by affine isometries, a stronger result was established by Shalom
in~\cite{Shalom}: One assumes only that the \emph{affine} $G$\hy
action does not almost have fixed points.
We replace Shalom's Hilbertian approach with an analogue of the geometric method used
in the splitting theorem of~\cite{Monod}.

(3)~This result can be reformulated in terms of the cohomology of
the associated linear $G$\hy representation $\ro$ on $B$ as
\[
    H^1(G,B)\ \cong\ \bigoplus_{i=1}^n H^1(G_i, B^{\ro(\prod_{j\neq i}G_j)}).
\]
It should be stressed that no such product formula holds in general.
Not only does it fail for more general Banach spaces (Example~\ref{E:product}),
but even for Hilbert space one needs at least Shalom's assumption mentioned above.
Compare the similar situation for the cohomological product formulas of~\cite{Shalom}
and~\cite{Burger-Monod}.
\end{rems}

\subsection{}
When $G$ is locally compact, we can as in the Lie case consider its lattices.
One then calls a lattice $\Gamma<G$ \emph{irreducible} if its projections to all $G_i$ are dense.
The above Theorem~\ref{thm_splitting} can be used to establish a superrigidity result for irreducible
lattices much in the way of~\cite{Shalom}. (The general idea to use irreducibility in order to
transfer results from $G_1\times \cdots \times G_n$ to $\Gamma$ was also illustrated
in~\cite{Burger-Mozes},
\cite{Burger-Monod},\cite{Monod-Shalom}; it seems to originate from the work
of Margulis and~\cite{BK}; lattices in products of completely
general locally compact groups were first studied by Shalom \cite{Shalom}.)
\begin{mthm}\label{C:super}
Let $\Gamma$ be an irreducible uniform lattice in a locally
compact $\sigma$\hy compact group $G=G_1\times \cdots \times G_n$.
Let $B$ be a superreflexive space with uniformly equicontinuous
affine $\Gamma$\hy action. Assume that the associated linear $\Gamma$\hy representation does not almost have invariant vectors.

Then there is a $\Gamma$\hy closed complemented affine subspace of
$B$ on which the $\Gamma$\hy action is a sum of actions extending continuously to $G$ and factoring through $G\to G_i$.
\end{mthm}
\begin{rem}\label{rem:C:super}
More precisely, the conclusion means that there are superreflexive
spaces $E_i$ endowed each with a continuous uniformly
equicontinuous affine $G$\hy action factoring through $G\to G_i$
and a $\Gamma$\hy equivariant affine continuous map
$\bigoplus_{i=1}^n E_i \to B$. Equivalently, the cocycle
$b:\Gamma\to B$ of the original $\Gamma$\hy action is cohomologous
to a sum $b_1+\cdots+ b_n$ of cocycles $b_i$ ranging in a subspace
$B_i\se B$ on which the \emph{linear} $\Gamma$\hy representation
extends continuously to a $G$\hy representation factoring through
$G_i$ and such that $b_i$ extends continuously to a cocycle $G\to
G_i\to B_i$ (with respect to the corresponding linear $G$\hy
representation). Moreover, $B_i\cong E_i$ as $G$\hy spaces.

If one disregards a component of $B$ where the linear $\Gamma$\hy representation ranges in a compact group of operators, this \emph{sum of actions} is actually just a direct sum $\bigoplus B_i\se B$ (see Remark~\ref{rem:compact:sum}).
\end{rem}

\begin{rem}
A \emph{uniform lattice} (in a locally compact group) is just a discrete cocompact subgroup; the theorem however also holds for certain non-uniform lattices, see Section~\ref{S:superrigidity} (Theorem~\ref{thm_super_prod}). Similar arguments allow us to generalise slightly Shalom's \emph{superrigidity for characters}, see Theorem~\ref{thm:characters}.
\end{rem}

\subsection*{Organization of the Paper}
In Section~\ref{S:prelims} we collect preliminary facts and lemmas
on uniformly convex/smooth and superreflexive Banach spaces,
linear representations and affine isometric on such spaces,
special properties of $L^{p}$\hy spaces,
and some general remarks and basic counter-examples.
In Section~\ref{S:DG} Theorem~\ref{T:T-FB-general} is proved.
Equivalence of properties~$(T)$ and~$(T_{L^{p}})$
(Theorem~\ref{T:TLp}) is proved in Section~\ref{S:TandTB}.
In Section~\ref{S:higher-rank} we discuss higher rank groups and
prove Theorem~\ref{T:higherrankFLp}. Section~\ref{sec_min} studies
minimal convex sets. Section~\ref{sec_product} addresses product groups
and proves the splitting theorem (Theorem~\ref{thm_splitting}); it also
proposes a proof of Theorem~\ref{T:higherrankFLp} that provides some evidence
for Conjecture~\ref{C:higherrankFB}. In Section~\ref{S:superrigidity}, we prove
Theorem~\ref{C:super}. Appendix~\ref{Appendix} describes Shalom's proof of
a generalized Howe--Moore theorem.

\subsection*{Acknowledgments}
We would like to thank D.~Fisher and G.A.~Margulis for their interest in this work and for letting us include
their argument for $2<p<2+\e$ in Section~\ref{S:FM}. We are indepted to A.~Naor for several helpful conversations
and to A.~Nevo for his remarks on the first manuscript. We are grateful to Y.~Shalom for letting us
give his proof of a ucus version of Howe--Moore (Theorem~\ref{T:HM}).

\section{Preliminaries} \label{S:prelims}
This section contains basic definitions, background facts and some
preliminary lemmas to be used in the proofs of our main results.

\subsection{Banach Spaces} \label{SS:Banach}

Let $B$ be a Banach space; unless otherwise specified, we take the reals as scalar field. We denote by $\sphere(B)=\Set{x\in B}{\|x\|=1}$ its unit sphere. For $x\in B$ and $r>0$ we denote by $\ball(x,r)$ and $\ol\ball(x,r)$ the open, respectively closed, ball of radius $r$ around $x$.

A Banach space $B$ is said to be \emph{strictly convex} if its unit sphere
does not contain straight segments, or equivalently if
$\|(x+y)/2\|<1$ whenever $x\neq y\in \sphere(B)$.
A Banach space $B$ is called \emph{uniformly convex} if
the \emph{convexity modulus} function
\begin{equation}\label{e:convexity-modulus}
    \delta(\gep)=\inf \Set{ 1-\|x+y\|/2 }{\|x\|,\|y\|\leq 1,\|x-y\|\ge\gep}
\end{equation}
is positive $\delta(\gep)>0$ whenever $\gep>0$.

We shall also use the notion of \emph{uniform smoothness} of Banach spaces,
which is easiest to define as the uniform convexity of the dual space $B^{*}$ (see~\cite[App.~A]{BL}).
Hence a Banach space $B$ is \emph{uniformly convex and uniformly smooth}
(hereafter abbreviated \emph{ucus}) if both $B$ and its dual $B^{*}$ are
uniformly convex.

\begin{facts} \label{F:ucus}
We refer to~\cite{BL} for the following:
\begin{enumerate}
\item
   The function $\gd (\gep )$ is non-decreasing and tends to $0$
   when $\gep$ tends to $0$. If $B$ is uniformly convex then
   $\gd(\gep )\to 0\Leftrightarrow\gep\to 0$.
\item
   Uniformly convex Banach spaces are reflexive.
   Hence the class of ucus Banach spaces is closed under taking duals.
   This class is also closed under the operations of taking
   closed subspaces and quotients.
\item
   If $B^{*}$ is strictly convex, in particular if $B$ is uniformly smooth,
   then every $x\in\sphere(B)$ has a unique supporting functional $x^{*}\in \sphere(B^{*})$,
   i.e. a unit functional with $\ip{x}{x^{*}}=1$.
\item
   If $B$ is ucus then the \emph{duality map} $*:\sphere(B)\to \sphere(B^{*})$, $x\mapsto x^{*}$,
   is a \emph{uniformly continuous}  homeomorphism with a uniformly continuous
   inverse.
\item \label{item:circumcentre}
  To any non empty bounded subset $E\se B$ of a reflexive strictly
  convex Banach space $B$, one can associate a unique point $C(E)\in B$,
  the \emph{circumcentre} of $E$ (a.k.a. the Chebyshev centre),
  defined as the unique $x\in B$ minimizing $\inf\{r>0: E\se\ol\ball(x,r)\}$.
\end{enumerate}
\end{facts}
The existence of $x=C(E)$ in~\eqref{item:circumcentre} follows
from weak compactness of closed bounded convex sets (i.e. from
reflexivity), whilst the uniqueness follows from uniform
convexity. Note that somewhat contrary to the intuition, it was
shown by V.~Klee~\cite{Klee60} that if $\dim (B)\geq 3$ and $B$ is
not a Hilbert space, then there exist a bounded subset $E\se B$
for which $C(E)$ does not belong to the closed convex hull of $E$.
The notion of circumcentre is also used in ${\rm CAT}(0)$
geometry.
For CAT$(0)$ spaces, the circumcentre $C(E)$ always lies in the closed convex hull of
$E$\footnote{Note that Hilbert spaces are,
in a sense, the most convex Banach spaces~-- they have the largest
possible modulus of continuity $\gd (\gep )$ among Banach spaces.
On the other hand, Hilbert spaces have the smallest possible
modulus of continuity among CAT$(0)$ spaces. Thus, in a sense,
CAT$(0)$ spaces are more convex then (non-Hilbertian) Banach
spaces}.

\smallskip

The following can be found e.g. in~\cite[A.6, A.8]{BL}:

\begin{thm}\label{F:superreflexive}
The following conditions on a topological vector space $V$
are equivalent:
\begin{enumerate}
\item $V$ is isomorphic to a uniformly convex Banach space.
\item $V$ is isomorphic to a uniformly smooth Banach space.
\item $V$ is isomorphic to a ucus Banach space.
\end{enumerate}
\end{thm}

The space $V$ is called \emph{superreflexive} if these equivalent condition hold. The class of superreflexive spaces is
closed under taking duals, closed subspaces and quotients of topological vector spaces.

\subsection{Linear Representations}

Let $V$ be a topological vector space.
We denote by $\GL(V)$ the group of invertible linear transformations
of $V$ which are continuous together with their inverses.

Following the standard terminology~\cite[Def.~2 of \S2 \no1]{BourbakiTG10},
a group $G$ of transformations of $V$ is \emph{uniformly equicontinuous}
(with respect to  the uniform structure deduced from the topological vector space structure)
if for any neighbourhood $U$ of $0\in V$ there exists a neighbourhood $W$ of $0$ such that
\begin{equation}\label{e:uniequi}
    x-y \in W \quad\Longrightarrow\quad \forall\, g\in G:\ g(x)-g(y)\in U.
\end{equation}
This definition will be applied to both linear groups, or affine groups.

\smallskip

For a topological vector space $V$, we denote by $N(V)$ the (a priori possibly empty)
set of norms on $V$ defining the given topology.
Elements of $N(V)$ will be called \emph{compatible} norms and are pairwise equivalent.

\medskip

The following key proposition is an equivariant version of Theorem~\ref{F:superreflexive}.
It enables us to reduce questions about uniformly
equicontinuous linear representations on superreflexive spaces to isometric
linear representations on ucus Banach spaces.
\begin{prop}[{\bf Invariant ucus norm}]\label{P:superreflex-ucus}
For a superreflexive topological vector space $V$
and a group of linear transformations $G$ of $V$, the following
conditions are equivalent:
\begin{enumerate}
\item
    $G$ is a uniformly equicontinuous group of linear transformations of $V$.
\item
    $G$ acts by uniformly bounded linear transformations
    with respect to any/all compatible norm on $V$.
\item
    $G$ acts by linear isometries with respect to some
    uniformly convex compatible norm on $V$.
\item
    $G$ acts by linear isometries with respect to some
    uniformly smooth compatible norm on $V$.
\item
    $G$ acts by linear isometries with respect to some
    uniformly convex and uniformly smooth compatible norm on $V$.
\end{enumerate}
\end{prop}
\begin{proof}
The main part of the proof is the implication ``[(3) and (4)]$\Rightarrow$(5)''; we begin by establishing this.

Let $N(V)$ denote the set of all compatible norms on $V$
equipped with the metric
\[
    d(\|\cdot\|_1, \|\cdot\|_2)=\sup_{x\neq 0}\Big| \log\frac{\|x\|_1}{\|x\|_2} \Big|.
\]
This is a complete metric space.
Let $N(V)^{G}$ stand for the closed subspace of $G$-invariant norms in $N(V)$.
Denoting by $\delta_{\|\cdot\|}$ the convexity modulus of $\|\cdot\|\in N(V)^{G}$,
the subset $N_{uc}(V)^{G}$ of uniformly convex
$G$-invariant norms on $V$ is given by the countable intersection
\[
    N_{uc}(V)^{G}=\bigcap_{n=1}^{\infty} O_{n},\qquad \text{where}\qquad
    O_{n}=\Set{\|\cdot\|\in N(V)^{G}}{\delta_{\|\cdot\|}(1/n)>0}.
\]
Observe that the sets $O_{n}$ are open.
If $\|\cdot\|_{0}$ is some fixed $G$-invariant compatible uniformly convex norm
(given in (3)) then any $\|\cdot\|\in N(V)^{G}$ can be viewed as a limit of uniformly convex
norms $\|\cdot\|+\e\|\cdot\|_0$ as $\e\searrow 0$.
Hence $N_{uc}(V)^{G}$ is a dense $G_{\delta}$ set in $N(V)^{G}$.

By duality between $N_{uc}(V^{*})^{G}$ and the set $N_{us}(V)^{G}$ of uniformly
smooth norms in $N(V)^{G}$, the latter is also a dense $G_{\delta}$ set in the Baire space
$N(V)^{G}$. In particular $N_{uc}(V)^{G}\cap N_{us}(V)^{G}$ is not empty, as claimed.

\smallskip

Now we observe that ``(1)$\Leftrightarrow$(2)'' follows from the definitions and that ``(5)$\Rightarrow$[(3) and (4)]''
as well as ``[(3) or (4) or (5)]$\Rightarrow$(2)'' are trivial. Moreover, proving ``(2)$\Rightarrow$(3)'' will also yield
``(2)$\Rightarrow$(4)'' by duality, using the fact that the dual to a superreflexive
space is superreflexive. Therefore it remains only to justify ``(2)$\Rightarrow$(3)'':

Let $\|\cdot\|$ be a compatible uniformly convex norm on $V$.
The corresponding operator norms $\|g\|=\sup_{x\neq 0} \|g x\|/\|x\|$
are uniformly bounded by some $C<\infty$.
Hence
\[
    \|x\|^{\prime}=\sup_{g\in G}\|g x\|
\]
defines a norm, equivalent to $\|\cdot\|$, and $G$-invariant.
It is also uniformly convex.
Indeed, if $\|x\|^{\prime}=\|y\|^{\prime}=1$ and $\|(x+y)/2\|^{\prime}>1-\alpha$
then for some $g\in G$
\[
    \|(gx+gy)/2\|>1-\alpha\qquad\text{whilst}\qquad
    \|gx\|\le \|x\|^{\prime}=1, \ \|gy\|\le\|y\|^{\prime}=1.
\]
Thus $\alpha\ge\delta_{\|\cdot\|}(\|gx-gy\|)\ge \delta_{\|\cdot\|}(\|x-y\|^{\prime}/C)$. Hence the convexity moduli
satisfy
\[
    \delta_{\|\cdot\|^{\prime}}(\e)\ge \delta_{\|\cdot\|}(\e/C)>0
    \qquad \text{for all}\quad \e>0.
\]
\end{proof}

If $G$ is a topological group, one should impose a continuity assumption
on linear $G$-representations on $V$, that is on homomorphisms $\ro:G\to\GL(V)$.
$\GL(V)$ is naturally equipped with the operator norm (which is
too strong for representation theory), and with the \emph{weak} and the \emph{strong} operator
topologies.
For uniformly equicontinuous representations the latter two topologies impose the same
continuity assumption:

\begin{lem}\label{L:cont-reps}
Let $G$ be a topological group, $V$ a superreflexive topological vector space, and $\ro:G\to\GL(V)$
a homomorphism. Then the following are equivalent.
\begin{enumerate}
\item $\ro$ is weakly continuous.
\item $\ro$ is strongly continuous.
\item The orbit maps $g\mapsto \ro(g)u$ are continuous.
\item The action map $G\times B\rightarrow B$ is jointly continuous.
\end{enumerate}
\end{lem}

Since there is an invariant complete norm on $V$, this is a special case of a well-known fact holding for all Banach spaces, see~\cite[3.3.4]{Monod_LNM} for references. We give an elementary proof in the present case.

\begin{proof}
Clearly it is enough to prove $(1)\Rightarrow(4)$.
Let $\|\cdot\|$ be a $\ro(G)$-invariant ucus norm on $V$.
Assume $g_n\rightarrow e\in G$ and $u_n\rightarrow u \in \sphere(B)$.
Then
\begin{eqnarray*}   |\ip{\ro(g_n)u_n}{u^*} -1| &\leq & \left| \ip{\ro(g_n)u_n}{u^*}-\ip{\ro(g_n)u}{u^*}\right|
        +\left|  \ip{\ro(g_n)u}{u^*}-1\right| \\
    &\le & \|u_n-u\| +|\ip{\ro(g_n)u}{u^*}-1| \rightarrow 0
\end{eqnarray*}
It follows that $\ro(g_n)u_{n}\rightarrow u$ because
\[
    \ip{\frac{\ro(g_{n})u+u}{2}}{u^{*}}\le \|\frac{\ro(g_{n})u+u}{2}\|
    \le 1-\delta(\|\ro(g_{n})u-u\|)
\]
and the left hand side tends to $1$.
\end{proof}

\subsection{Invariant complements}
One of the convenient properties of Hilbert spaces is the
existence of a canonical complement $M^{\perp}$ to any closed
subspace $M$. Recall that a closed subspace $X$ of a Banach space
$B$ is called complemented if there is another closed subspace
$Y\leq B$ such that $B=X\oplus Y$ algebraically and topologically.
This is equivalent to each of the following:
\begin{itemize}
\item
  There is a continuous linear projection from $B$ to $X$.
\item
  There is a closed subspace $Y$ and a continuous linear
  projection $p:B\to Y$ with $\ker (p)=X$.
\end{itemize}
A classical theorem of Lindenstrauss and Tzafriri says that every
infinite dimensional Banach space which is not isomorphic to
a Hilbert space, admits a non-complemented closed subspace~\cite{Lindenstrauss-Tzafriri_compl}.
However, for any uniformly equicontinuous linear representation $\ro$ of
a group $G$ on a superreflexive space $B$, the
subspace of invariant vectors $B^{\ro(G)}$ admits a canonical
complement, described below.

In view of Proposition~\ref{P:superreflex-ucus} we may assume that
the representation is linear isometric
with respect to a ucus norm on $B$, which allows
to use the duality map of the unit spheres $*:\sphere(B)\to \sphere(B^{*})$.

Given any linear representation $\ro:G\to\GL(V)$ there is an associated dual
(or contragradient) linear $G$-representation $\ro^{*}:G\to\GL(V^{*})$
defined by
\[
    \ip{x}{\ro^{*}(g)y}=\ip{\ro(g^{-1})x}{y}\qquad (g\in G,\ x\in V,\ y\in V^{*}).
\]
If $B$ is a Banach space and $\ro:G\to\Orth(B)$ is a linear
isometric representation, then so is its dual
$\ro^{*}:G\to\Orth(B^{*})$, where $B^{*}$ is equipped with the
dual norm. Hence the dual to a uniformly equicontinuous
representation on a superreflexive space is also of the same type.
\begin{ob} \label{Ob:fixed-points}
If $B$ is a ucus Banach space and $\ro:G\to\Orth(B)$, then the duality map
$*:\sphere(B)\to \sphere(B^{*})$ between the unit spheres
intertwines the actions of $\ro(G)$ and $\ro^*(G)$.
In particular it maps the set of $\ro(G)$\hy fixed unit vectors to the set of
$\ro^*(G)$\hy fixed unit vectors.
\end{ob}

\begin{prop}\label{P:linear-splitting}  
Let $\ro$ be a uniformly equicontinuous linear representation of $G$
on a superreflexive space $B$, let $B^{\ro(G)}$ denote
the subspace of $\ro(G)$\hy fixed vectors in $B$, and let
$B'=B'(\ro)$ be the annihilator of $(B^*)^{\ro^*(G)}$ in $B$.
Then
\[
    B=B^{\ro(G)}\oplus B'(\ro).
\]
Furthermore, the decomposition is canonical in the
following sense: If we denote by $p(\ro)$ and $p'(\ro)$ the associated
projections, then for every morphism of uniformly equicontinuous linear representations
$\phi:(B_1,\ro_1)\rightarrow (B_2,\ro_2)$, the following diagrams are commutative:
\begin{equation}\label{functoriality}
     \begin{CD}
      B_1  @>{\phi}>>  B_{2}  \\
       @V{p(\ro_{1})}VV   @V{p(\ro_{2})}VV\\
        B_1 @>{\phi}>>  B_2
     \end{CD}
    \qquad\qquad
     \begin{CD}
      B_1  @>{\phi}>>  B_{2}  \\
       @VV{p'(\ro_{1})}V   @VV{p'(\ro_{2})}V\\
        B_1 @>{\phi}>>  B_2
     \end{CD}
\end{equation}
\end{prop}

\begin{rem}\label{rem:no_complement}
The conclusion fails if we drop the superreflexivity assumption, see Example~\ref{exo:no_complement}.
\end{rem}

\begin{proof}[Proof of the proposition]
Choose a $G$-invariant uniformly convex and uniformly smooth norm on $B$,
and the dual one on $B^{*}$ (Proposition~\ref{P:superreflex-ucus}).
For any unit vector $x\in B^{\ro(G)}$ and arbitrary $y\in B^{\prime}$
\[
  1=\ip{x}{x^{*}}=\ip{x-y}{x^{*}}   \le \|x-y\|\cdot\|x^{*}\|=\|x-y\|.
\]
Thus $B^{\ro(G)}\cap B^{\prime}=\{0\}$ and $B^{\ro(G)}\oplus
B^{\prime}$ is a closed subspace in $B$. It is also dense in $B$.
Indeed if $\lambda\in B^{*}$ is a unit vector vanishing on
$B^{\prime}$ it cannot vanish on $B^{\ro(G)}$, because $\lambda\in
(B^*)^{\ro^*(G)}$ by the Hahn--Banach theorem, and hence
$\lambda^*\in B^{\ro(G)}$ and $\ip{\lambda^{*}}{\lambda}=1$. Thus
$B^{\ro(G)}\oplus B^{\prime}=B$.

The last assertion follows from the fact that $\phi(B_1^{\ro_1}) \se
B_2^{\ro_2}$, and  $\phi^*((B_2^*)^{\ro_2})\se (B_1^*)^{\ro_1}$ yields
$\phi(B_2')\se B_1'$.
\end{proof}

\begin{cor}\label{normelizers}
The decomposition $B=B^{\ro(G)}\oplus B^{\prime}$ is preserved
by the normalizer of $\ro(G)$ in $\GL(B)$.
\end{cor}

\begin{cor}\label{C:split:ucus}
Let $G=G_1\times  G_2$ be any product of two groups and $B$ a superreflexive space
with a uniformly equicontinuous  linear $G$\hy representation $\ro$.
Then there is a canonical $G$\hy invariant decomposition
\[
    B = B^{\ro(G)} \oplus B_0 \oplus B_1 \oplus B_2
\]
such that $B^{\ro(G_i)} = B^{\ro(G)} \oplus B_i$ for $i=1,2$.
\end{cor}

\begin{prop}\label{stam}
Let $\ro$ be a uniformly equicontinuous linear $G$-representation on a
superreflexive space $B$. Then
\begin{enumerate}
\item
    $B^{\ro (G)}$ is isomorphic to $B/B'$ as topological vector spaces.
\item
    $B'$ is isomorphic to $B/B^{\ro (G)}$ as $G$\hy representations.
\item
    $(B^{\ro(G)})^*$ is isomorphic to $B^*/(B^*)'$ as topological vector spaces.
\item
    $(B')^*$ is isomorphic to $(B^*)'$ as $G$\hy representations.
\item
    $B'$ almost has invariants if and only if $(B^*)'$ almost has invariants.
\item
    If $\ 0\rightarrow A\rightarrow B\rightarrow C\rightarrow 0\ $ is an
    exact sequence of uniformly equicontinuous linear $G$\hy representations
    on superreflexive spaces, then $B'$ almost has invariant vectors
    if and only if $A'$ or $C'$ do.
\end{enumerate}
If $B$ is equipped with a compatible uniformly convex and uniformly smooth
$G$\hy invariant norm, then the natural isomorphisms in {\rm (1)} and {\rm (3)}
are isometric.
\end{prop}

\begin{proof}
Equip $B$ with a $G$-invariant ucus norm (Proposition~\ref{P:superreflex-ucus}).

By the open mapping theorem the maps $p:B\to B^{\ro(G)}$ and $p':B\to B^{\prime}$
induce isomorphisms of topological vector spaces
\[
    (1)\quad \tilde{p}:B/B'\rightarrow B^{\ro(G)},\qquad
    (2)\quad \tilde{p'}:B/B^{\ro(G)}\rightarrow B'.
\]
By~\ref{Ob:fixed-points} $(B^{\ro (G)})^*$ is $(B^{*})^{\ro^{*}(G)}$ and the latter is
isomorphic to  $B^*/(B^*)'$. This proves (3).

To see that (1) and (3) are isometric (with respect to the norms corresponding to
any ucus $G$-invariant norm on $B$) we note that the isomorphisms above satisfy
$\|\tilde{p}^{-1}\|,\|\tilde{p^{\prime}}^{-1}\|\le 1$
by the definition of the norm on a quotient space.
Furthermore, for $v \in \sphere(B^{\ro (G)})$, we have $v^*\in \sphere((B^*)^{\ro^{*}(G)})$,
hence
\begin{eqnarray*}
    \| \tilde{p}^{-1}(v) \|_{B/B'}&=&\inf \Set{\|v+v'\|_{B} }{v'\in B'}\\
        &\geq& \inf\Set{\ip{v+v'}{v^{*}}}{v'\in B'} = \ip{v}{v^{*}}=1.
\end{eqnarray*}
Hence $\tilde{p}$ is an isometry $B/B'\cong B^{\ro(G)}$.
Similarly, $(B^{\ro (G)})^*=(B^{*})^{\ro^{*}(G)}\cong B^*/(B^*)'$.

In general Banach spaces the dual $E^{*}$ of a subspace $E<F$ is isometric
to the quotient $F^{*}/E^{\perp}$ by the annihilator $E^{\perp}<F^{*}$ of $E$.
Thus with respect to a ucus norm on $B$ and the above spaces,
$(B')^*$ is isometric to $B^*/(B^*)^{\ro^{*}(G)}$ as Banach spaces, while
the latter is isomorphic to $(B^*)'$ as a topological vector space by (2).
Whence (4).

(5) Assume that there exist $x_{n}\in \sphere(B^{\prime})$ with $\diam(\ro(K)\cdot x_{n})\to0$.
The uniformly continuous map $*:\sphere(B)\to\sphere(B^{*})$
takes vectors $x_{n}\in\sphere(B^{\prime})$ to vectors $x^{*}_{n}\in \sphere(B^{*})$
with $\diam(\ro^{*}(K)\cdot x^{*}_{n})\to0$.
Since $\{x^{*}_{n}\}$ are uniformly separated from $(B^{*})^{\ro^{*}(G)}$
their  normalized projection $y^{*}_{n}$ to $(B^{*})^{\prime}$ still satisfy
$\diam(\ro^{*}(K)\cdot y^{*}_{n})\to0$.

(6) As $A'$ maps into $B'$, if $A'$ almost has invariants, then so
does $B'$. If $C'$ almost has invariants, then so does $(C^*)'$,
hence $(B^*)'$, hence $B'$. On the other hand, assume $B'$
almost has invariant unit vectors $b_n$. Assume for simplicity that $A=A'$, $B=B'$
and $C=C'$. Note that $C$ is isomorphic to $B/A$,
and denote by $\pi:B\rightarrow C$ the projection. Then either $\pi(b_n)$ converges
to $0\in C$, then there exist $a_n$ such that $b_n-a_n$ converges
to $0\in B$, and the normalized sequence $(\frac{a_n}{\| a_n\|})$ is almost
invariant in $A$, or there exist a subsequence $b_{n_k}$ with
$\inf_k \|\pi(b_{n_k})\| > 0$, and then the normalized sequence
$(\frac{\pi(b_{n_k})}{\|\pi(b_{n_k})\|})$ is almost invariant in $C$.
\end{proof}

\begin{rem} \label{R:TB}
For ucus Banach space $B$, Definition~\ref{D:TB} of
property~$(T_{B})$ can be rephrased as follows: For any
representation $\ro:G\to\Orth(B)$, the restriction
$\ro^{\prime}:G\to\Orth(B^{\prime})$ of $\ro$ to the invariant
subspace $B^{\prime}$ complement to $B^{\ro(G)}$ does not almost
have invariant vectors, i.e. for some compact $K\se G$ and
$\gep>0$
\[
    \forall\, v\in\sphere(B^{\prime})\quad\exists\, g\in K\ \text{s.t.}\
    \|\ro(g)v-v\|\ge\gep.
\]
\end{rem}
Hence item (4) gives:
\begin{cor}\label{C:TB<=>TB}
Let $B$ be a ucus Banach space, and $G$ be a locally compact group.
Then $G$ has property~$(T_{B})$ iff it has~$(T_{B^{*}})$.
\end{cor}

\subsection{Affine Actions} \label{SS:affisoms}
The affine group $\Aff(V)$ of a real affine space $V$ (a vector space who forgot its origin)
consists of invertible maps satisfying:
\[
    T(t \cdot x + (1-t)\cdot y)=t\cdot T(x)+(1-t)\cdot T(y),
    \qquad(t\in\BR,\ x,y\in V)
\]
The group $\Aff(V)$ is a semi-direct product $\Aff(V)=\GL(V)\ltimes V$,
i.e. an invertible affine map $T$ has the form $T(x)=Lx+b$ where
$L\in \GL(V)$ is linear invertible.

An affine action of a group $G$ on $V$, i.e. a homomorphism $G\to\Aff(V)$, has the form
\[
    g\cdot x= \ro(g)x + c(g),
\]
where $\ro:G\to\GL(V)$ is a linear $G$-representation (we
call it the \emph{linear part} of the action) and  $c:G\to B$ is a
$\ro$\hy\emph{cocycle}, namely an element of the Abelian group
\begin{equation}\label{e:ro-cocycles}
    Z^{1}(\ro)=\Set{c:G\to V}{c(gh)=\ro(g)c(h)+c(g),\ \forall\, g,h\in G}.
\end{equation}
The group $Z^{1}(\ro)$ of $\ro$\hy cocycles contains the subgroup of
$\ro$\hy\emph{coboundaries}
\begin{equation}\label{e:ro-coboundaries}
    B^{1}(\ro)=\Set{ c(g)=v-\ro(g)v}{v\in V}.
\end{equation}
$Z^{1}(\ro)$ describes all affine $G$\hy actions on $V$ with linear part $\ro$,
and $B^{1}(\ro)$ corresponds to those affine actions which have a $G$\hy fixed
point (namely $v$ in~\eqref{e:ro-coboundaries}).
This description involves the choice of reference point~-- the origin~-- in the space.
Two cocycles differing by a coboundary can be though of defining the same affine
action viewed from different reference points.
The first cohomology of $G$ with $\ro$\hy coefficients is  the Abelian group
\[
    H^{1}(G,\ro)=Z^{1}(\ro)/B^{1}(\ro).
\]
It describes different types of actions in the above sense.
$H^{1}(G,\ro)=0$ iff any affine $G$\hy action on $V$ with linear part $\ro$
has a fixed point.

\medskip

For a Banach space $B$ denote by $\Isom(B)$ the group of isometries of
$B$ as a metric space.
It is a classical theorem of Mazur--Ulam that any surjective
isometry $T$ of a (real) Banach space $B$ is necessarily affine
$T(x)=L x +c$ with linear part $L\in \Orth(B)$ being isometric.
(This theorem is elementary when $B$ is strictly convex; comopare
Lemma~\ref{iso-affine}).
Hence $\Isom(B)=\Orth(B)\ltimes B$.

\medskip

Now suppose that $V$ is a superreflexive topological vector space.
Recall that a group $G$ of affine self maps is uniformly equicontinuous
if it satisfies~\eqref{e:uniequi}. This condition is equivalent to uniform equicontinuity of the
linear part $\ro:G\to\GL(V)$.
\begin{prop}\label{P:superreflex-affine}
For a superreflexive topological vector space $V$
and a group of transformations $G$ of $V$ the following
conditions are equivalent:
\begin{enumerate}
\item
    $G$ is uniformly equicontinuous group of affine transformations of $V$.
\item
    $G$ acts by uniformly Lipschitz affine transformations
    with respect to any/all compatible norms on $V$.
\item
    $G$ acts by affine isometries with respect to some
    compatible norm on $V$.
\item
    $G$ acts by affine isometries with respect to some
    uniformly convex and uniformly smooth compatible norm on $V$.
\end{enumerate}
\end{prop}

\begin{proof}
Apply Proposition~\ref{P:superreflex-ucus}
to the linear part of the affine action, using Mazur--Ulam to deduce in (3) that the action
is affine.
\end{proof}

\medskip

If $G$ is a topological group acting by affine transformations on a topological vector space $V$,
continuity of the action
\[
    G\times V\to V,\qquad g\cdot x=\ro(g)x+c(g)
\]
is equivalent to continuity of the linear part $G\times V\to V$ and the continuity
of the cocycle $c:G\to V$.
Indeed $c(g)=g\cdot 0$, and $\ro(g)x=g\cdot x-c(g)$.

Hence in the context of topological groups, affine actions should be assumed
continuous, and $Z^{1}(G,\ro)$ will include only continuous cocycles $c:G\to V$
(we assume that the linear part $\ro$ is continuous as well).
If $G$ is a locally compact $\sigma$\hy compact group, then $Z^{1}(\ro)$ has a
natural structure of a Fr\'echet space with respect to the family of semi-norms
\[
    \|c\|_{K}=\sup_{g\in K} \|c(g)\|_V
\]
where $K\se G$ runs over a countable family of compact subsets
which cover $G$ and $\|\cdot\|_V$ is a norm inducing the topology of $V$. Moreover, if $G$ is compactly generated (e.g. if
$G$ has property~$(T)$) say by $K_{0}$, then $\|c\|_{K_{0}}$ is a norm on
$Z^{1}(\ro)$ (note that any cocycle $c\in Z^{1}(\ro)$ is
completely determined by its values on a generating set), and
$Z^{1}(\ro)$ is a Banach space with respect to this norm.
We remark that in general $B^{1}(\ro)$ is not closed in $Z^{1}(\ro)$
(this is the idea behind the~$(F_{B})\Rightarrow (T_{B})$ argument
of Guichardet~-- see Section~\ref{S:DG}).

\medskip

\begin{lem}\label{L:bounded}
For a uniform equicontinuous affine action of a group $G$ on a
superreflexive space $B$, the following are equivalent:
\begin{enumerate}
\item There exists a bounded $G$\hy orbit.
\item All $G$\hy orbits are bounded
\item $G$ fixes a point in $B$.
\item $G$ preserves a (Borel regular) probability measure on $B$.
\end{enumerate}
\end{lem}
Note that the notion of a subset $E\se V$ being bounded, means that
for any open neighbourhood $U$ of $0\in V$ there is some $t\in\BR$
so that $E\se t\cdot U$.
This notion agrees with the notion of being bounded with respect to any
compatible norm on $V$.
\begin{proof}
Introduce a $G$-invariant uniformly convex norm on $V$ (Proposition~\ref{P:superreflex-affine}).
The only non-trivial implications are ${\rm (4)}\Rightarrow{\rm
(1)}\Rightarrow{\rm (3)}$. For the first, let $\mu$ be a $G$\hy
invariant probability on $B$. Since $B$ is a countable union of
closed bounded sets, there is a closed bounded set $A\se B$ with
$\mu(A)>1/2$. For all $g\in G$ we have $\mu(gA)>1/2$ hence $gA
\cap A\neq \varnothing$. It follows that the $G$\hy orbit of every
point of $A$ is bounded.

The latter implication follows by considering the circumcentre (compare Section~\ref{SS:Banach})
 of the given bounded $G$\hy orbit.
\end{proof}

\begin{prop}\label{P:generalities}
Let $B$ be a ucus Banach space. Then
\begin{enumerate}
\item
    Any finite (or compact) group has properties~$(T_{B})$ and~$(F_{B})$.
\item
    Properties~$(T_{B})$ and~$(F_{B})$ pass to quotient groups.
\item
    If $G=G_{1}\times\dots\times G_{n}$ is a finite product of topological groups
    then $G$ has property~$(T_{B})$ (resp.~$(F_{B})$) iff
    all $G_{i}$ have this property.
\end{enumerate}
\end{prop}

\begin{proof}
(1) and (2) are straightforward, (3) follows from Corollary~\ref{C:split:ucus}.
\end{proof}

\subsection{Special Properties of $L^{p}(\mu)$\hy Spaces}\label{SS:Lp}
In this section we collect some special properties of the Banach spaces
$L^{p}(\mu)$ which will be used in the proofs.

By an $L^{p}(\mu)$, or $L^{p}(X,\mu)$ space we mean the usual space of equivalence
classes (modulo null sets) of measurable $p$\hy integrable functions $f:X\to \BR$,
where $\mu$ is a positive $\sigma$\hy finite measure defined on
a standard Borel space $(X,\mathcal{B})$.
If $1<p<\infty$ then $L^{p}(\mu)$ is ucus, whilst
$L^{1}(\mu)$ and $L^{\infty}(\mu)$ are not (they are not even strictly convex).
For $1\le p<\infty$ the dual to $L^{p}(\mu)$ is $L^{q}(\mu)$ where $1<q\le\infty$ is determined
by $q=p/(p-1)$.

The space $L^{p}([0,1], \text{Lebesgue})$ is usually denoted by $L^{p}$.
Any $L^{p}(\mu)$\hy space with \emph{non-atomic} finite or $\sigma$\hy finite measure
$\mu$ is isometrically isomorphic to $L^{p}$.
Indeed let  $\varphi\in L^{1}(\mu)$ be a strictly positive measurable
function with integral one and let $\mu_{1}$ be given by $d\mu_{1}=\varphi\,d\mu$.
Then
\[
    f\in L^{p}(\mu) \mapsto f\cdot \varphi^{-1/p} \in L^{p}(\mu_{1})
\]
is a surjective isometry. Since any non-atomic standard
probability spaces is isomorphic to $[0,1]$ as a measure space,
$L^{p}(\mu_{1})\cong L^{p}$. If $\mu$ is purely atomic then a
similar argument gives an isomorphism of $L^{p}(\mu)$ with finite
or infinite dimensional $\ell^{p}$ space. A general $L^{p}(\mu)$
space is therefore isometrically isomorphic to a direct sum of
$L^{p}$ and $\ell^{p}$ components.

More generally, for another Banach space $B$, one defines the spaces
$L^{p}(\mu,B)$ of $B$\hy valued function classes by means of the Bochner integral.
We refer the reader to~\cite{Diestel-Uhl} for details; we recall here that the dual
of $L^{p}(\mu,B)$ is $L^q(\mu, B^*)$ through the natural pairing for all $1\leq p < \infty$,
but only when $B$ has the Radon--Nikod\'ym property~-- this includes all ucus spaces (see again~\cite{Diestel-Uhl}).
These spaces will be used in Section~\ref{sec_induction} in order to \emph{induce} isometric
(linear or affine) actions.

\medskip

Banach~\cite{Banach} and Lamperti~\cite{Lamperti} (see
also~\cite[Theorem~3.25]{FJ}) classified the linear isometries of
$L^{p}(\mu )$ as follows.
\begin{thm}[Banach, Lamperti]  \label{T:Banach-Lamperti}
For any $1<p<\infty$ where $p\neq 2$, any linear isometry $U$
of $L^{p}(X,\mathcal{B},\mu )$ has the form
\[
    Uf(x)=f(T(x))h(x)\left(\frac{dT_*\mu}{d\mu}(x)\right)^{\frac{1}{p}}
\]
where $T$ is a measurable, measure class preserving map of $(X,\mu)$,
and $h$ is a measurable function with $|h(x)|=1$ a.e.
\end{thm}
Let $\mu=\mu_{a}+\mu_{c}$ be the decomposition of $\mu$ into
its atomic and continuous parts ($\mu_{a}=\mu|_{A}$ where $A\se X$
is the (at most countable) set of atoms of $\mu$).
Then
\[
    L^{p}(\mu)=L^{p}(\mu_{c})\oplus L^{p}(\mu_{a})\cong L^{p}\oplus \ell^{p}(A)
    \quad\text{or just}\quad \ell^{p}(A),
\]
the latter case occurs if $\mu=\mu_{a}$ is a purely atomic measure.
Note that it follows from Banach--Lamperti theorem that this decomposition
is preserved by any linear isometry of $L^{p}(\mu)$.
As $\ell^{p}(A)$ has a much smaller group of linear (or affine) isometries
than $L^{p}$ we could restrict our attention only to the latter.
However we shall not make use of this ``simplification``.

Another useful tool in the study of $L^{p}$\hy spaces is the Mazur map.
\begin{thm}[ {\cite[Theorem~9.1]{BL}} ] \label{T:Mazur}
Let $\mu$ be a $\sigma$\hy finite measure. For any $1\leq p,q<\infty$
the \emph{Mazur map} $M_{p,q}:L^{p} (\mu )\to L^{q} (\mu)$ defined by
\[
    M_{p,q}(f)={\rm sign}(f)\cdot |f|^{\frac{p}{q}}
\]
is a (non-linear) map which induces a uniformly continuous homeomorphism
between the unit spheres $M_{p,q}:\sphere(L^{p}(\mu ))\to \sphere(L^{q}(\mu ))$.
\end{thm}
(Note that if $p,q\neq 1$ and $p^{-1}+q^{-1}=1$ then the restriction of $M_{p,q}$ to the
unit spheres is just the duality map $*:\sphere(L^{p}(\mu))\to \sphere(L^{p}(\mu)^{*})$).

In the proofs of Theorems~\ref{T:TLp} and~\ref{T:higherrankFLp},
the results for subspaces and quotients are deduced from the
$L^{p}(\mu)$ case using the following theorem of Hardin about
extension of isometries defined on subspaces of $L^{p}(\mu)$. The
formulation we give here is not quite identical to the original,
but it easily follows from it and from its proof
(see~\cite[Theorem~4.2]{Hardin} or~\cite[Theorem~3.3.14]{FJ}).
\begin{thm}[Hardin]
Let $(X,\mathcal{B},\mu)$ be a measure space. For every closed
subspace $F\se L^{p}(X,\mu)$, there is a canonical extension
$F\se\ti F\se L^{p}(\mu)$ which is isometric to
$L^{p}(X',\mu')$ for some other measure space $(X',\mu')$.
Furthermore, if $1<p\notin 2\BZ$, then every linear isometry
$U:F\rightarrow L^{p}(Y,\nu)$ extends uniquely to a surjective
linear isometry $\ti U:\ti F\rightarrow \widetilde{UF}\subseteq
L^{p}(Y,\nu)$.
\end{thm}
\begin{rem}
If $\mathcal{B}'\leq\mathcal{B}$ is the minimal sub $\gs$\hy algebra
with respect to which all the functions in $F$ are measurable,
then $\ti F=L^{p}(X,\mathcal{B}^{\prime},\mu )$.
\end{rem}
A straightforward consequence is the following:
\begin{cor}\label{C:Hardin}%
Let $1 < p\notin 2\BZ$, and let $F\se  L^p(X,\mu)$ be a closed
subspace. Let $\ro$ be a linear isometric representation of the
group $G$ on $F$. Then there is some linear isometric $G$\hy
representation $\ro'$ of $G$ on some other space $L^{p}(X',\mu')$,
and a linear $G$\hy equivariant isometric embedding
$F\hookrightarrow L^{p}(X',\mu')$.
\end{cor}
Another important fact about $L^{p}(\mu)$\hy spaces, this time for
$p\in (0,2]$, is that $B=L^{p}(\mu)$ has an embedding
$j:B\to\cH$ into the unit sphere of a Hilbert space so
that $\langle j(x),j(y)\rangle =\|x-y\|^{p}$. Having such an
embedding is equivalent (via the classical result of I.J.~Schoenberg,
see~\cite{BHV}) to the following:
\begin{prop} \label{F:pos-def}
For $0<p<2$ and any $s>0$ the function $f\mapsto e^{-s\|f\|^{p}}$
is positive definite on $L^{p}(\mu)$, i.e. for any finite collection
$f_{i}\in L^{p}(\mu)$ and any $\lambda_{i}\in \BC$
\[
    \sum_{i,j} e^{-s \|f_{i}-f_{j}\|^{p}}\lambda_{i}\overline{\lambda_{j}} \ge 0.
\]
\end{prop}

In fact, more is known: It was shown by Bretagnolle, Dacunha-Castelle
and Krivine~\cite{BDCK} (c.f.~\cite[5.1]{WW}) that, for $1\leq p\leq 2$, a
Banach space $X$ is isometric to a closed subspace of $L^{p}(\mu )$
iff $e^{-s\|\cdot\|^p}$ is a positive definite function on $X$ for any $s>0$.

\subsection{Some Easy Counterexamples and Remarks}
\label{SS:couterexamples}
%
\begin{exam}[$(T_{B})\not\Rightarrow (F_{B})$]  \label{T+notF}
Let $B$ be a Banach space with $\Orth(B)\cong \BZ/2\BZ$, i.e.
a space where the only linear isometries are the identity and
the antipodal map $x\mapsto -x$.
A trivial example of such a space is the line $B=\BR$,
but it is not hard to construct such spaces of arbitrary dimensions
even within the class of ucus Banach spaces (by considering e.g. sufficiently asymetric
convex sets in Hilbert space and choosing the corresponding norm).
Clearly for such a space  any group has property~$(T_B)$.
However the groups $\BZ$ or $\BR$ or any group $G$
with sufficiently large Abelianization $G/[G,G]$ would fail to have property
$(F_B)$ for it would admit an isometric action by translations:
$n\cdot x:=x+nx_0$
where $0\neq x_0\in B$ is arbitrary.
However groups with trivial Abelianization would also have
property~$(F_{B})$ on such an asymmetric Banach space $B$.
\end{exam}
\begin{exam}[$(T)\not\Rightarrow (F_{B})$]\label{E:L10}
Suppose $G$ is locally compact non-compact (e.g. $G=\SL_3(\BR )$ or $G=\SL_3(\BZ)$).
Fix a Haar measure on $G$ and let
\[
    B=L^{1}_{0}(G)=\Set{ f\in L^{1}(G)}{\int f\,dg=0}
\]
be the codimension one subspace of functions with $0$ mean. Then
$B$ is isometric to the affine subspace $\Set{ f\in L^{1}(G)}{\int
f=1}$ on which $G$ acts isometrically by translations without
fixed points. Hence $G$ does not have property~$(F_B)$. This
Banach space is not ucus. Notice that in this example all orbits
are bounded regardless of $G$.
\end{exam}

\begin{rem}
Haagerup and Przybyszewska \cite{Haagerup-Przybyszewska} showed
that any locally compact group $G$ admits a proper isometric
action on the strictly convex space $\bigoplus_{n=1}^\infty
L^{2n}(G)$.
\end{rem}

\begin{exam}[$(T)\not\Rightarrow (T_{B})$]\label{2.21}
Let $G$ be as in Example~\ref{E:L10}. Consider the space
$B=C_0(G)$ of continuous real valued functions on $G$ which tend
to $0$ at $\infty$ with the sup ($L^\infty$\hy ) norm. The action
of $G$ on $B$ by translations is a linear isometric action. A
function $f\in B$ which decays very slowly forms an ``almost
invariant vector''. On the other hand there are no non-zero
invariant vectors. Hence $G$ does not have property~$(T_B)$.
\end{exam}
\begin{rem}
Since any separable Banach space is a quotient of $\ell^{1}$,
Example~\ref{2.21} shows that case (iii) of Theorem~\ref{T:TLp}
cannot be extended to $p=1$.
\end{rem}
\begin{exam}[Remarks~\ref{R:thm_splitting}(2)]\label{E:product}
Let $G=G_1\times G_2$ be any product of non-compact locally
compact groups (e.g. $G=\BZ\times \BZ$). Let $B=L_0^1(G)$ as in
Example~\ref{E:L10}. Then $H^1(G,B)\neq 0$, but there are no
non-zero $G_i$\hy fixed vectors in the associated linear
representation. Thus the product formula of
Remarks~\ref{R:thm_splitting}(2) cannot hold for $B$.
\end{exam}

Let us make some remarks about Kazhdan's property~$(T)$ and
property~$(T_{B})$ as in~\ref{D:TB} and~\ref{R:TB}. Given a
unitary representation $(\ro,\cH)$ of a locally compact
group $G$, a compact subset $K\se G$ and $\gep>0$, one says
that a vector $0\neq v\in\cH$ is $(K,\gep)$\hy \emph{
almost invariant} if
\[
    \sup_{g\in K} \|\ro(g)v-v\|<\gep\cdot\|v\|.
\]
A locally compact group $G$ has Kazhdan's property~$(T)$
if and only if it satisfies the following equivalent conditions:
\begin{enumerate}
\item
    For any unitary $G$\hy representation $(\ro, \cH)$
    there exists a compact $K\se G$ and an $\gep>0$
    so that the $G$\hy representation $\ro^{\prime}$ on
    $(\cH^{\ro(G)})^{\perp}\cong \cH/\cH^{\ro(G)}$
    has no $(K,\gep)$\hy almost invariant vectors.
\item
    There exist a compact $K\se G$ and an $\gep>0$ so that all non-trivial irreducible unitary $G$\hy representations
    $(\ro, \cH)$ have no $(K,\gep)$\hy almost invariant vectors.
\item
    There exist a compact $K\se G$ and an $\gep>0$ so that for all unitary $G$\hy representations $(\ro, \cH)$
    the $G$\hy representation $\ro^{\prime}$ on
    $(\cH^{\ro(G)})^{\perp}\cong \cH/\cH^{\ro(G)}$
    has no $(K,\gep)$\hy almost invariant vectors.
\end{enumerate}
In the above, (3) clearly implies both (1) and (2). In showing
(1)$\Rightarrow$(3) one uses the fact that the category of Hilbert
spaces and unitary representations is closed under $\ell^{2}$ sums
and $L^{2}$\hy integration. The fact that any unitary
representation decomposes as an $L^{2}$\hy integral of irreducible
ones gives (2)$\Rightarrow$(3).

\begin{rem}\label{R:TBuniform}
Definition~\ref{D:TB} (Remark~\ref{R:TB}) of property~$(T_{B})$ is
modeled on (1) above. There does not seem to be any reasonable
theory of irreducible representations (and decomposition into
irreducibles) for Banach spaces other than Hilbert ones. Hence
form (2) of property~$(T)$ does not seem to have a Banach space
generalization. As for (3), for any given $1<p<\infty$ the class
of $L^{p}(\mu)$\hy spaces is closed under taking $\ell^{p}$\hy
sums (and $L^{p}$\hy integrals) and hence for groups with
property~$(T_{L^{p}})$ an analogue of (3) holds, namely there
exist $K\se G$ and $\gep>0$ which are good for all
$\ro:G\to\Orth(L^{p})$.
Also, if a group $G$ has property~$(T_{B})$ for \emph{all} ucus Banach spaces $B$
(conjecturally all higher rank groups and their lattices) then for every ucus
Banach space $B$ there is $(K,\gep)$ which is good for all linear isometric representations
$G\to\Orth(B)$.
This uses the fact that $L^{2}(\mu,B)$ is a ucus
if $B$ is (see Lemma~\ref{lem_ind_uc} below).
\end{rem}

Finally, we justify Remark~\ref{rem:no_complement}:

\begin{exam}\label{exo:no_complement}
Let $G$ be a discrete group and consider the Banach space
$B=\ell^\infty(G)$ with the (linear isometric) regular $G$\hy
representation $\ro$. Then one shows that the space $B^{\ro(G)}$
(which consists of the constant functions) admits a $G$\hy
invariant complement (if and) \emph{only if $G$ is amenable}.
Indeed, the Riesz space (or Banach lattice) structure of $B$
allows to take the ``absolute value'' of any linear functional on
$B$; renormalizing the absolute value of any non-zero invariant
functional would yield an invariant mean on $G$. Alternatively,
one can argue similarly on the Banach space of continuous
functions on any compact topological $G$\hy space.

We point out that nevertheless the space $B'$ is well-defined for
any topological vector $G$\hy space $B$; in the case at hand, we
have $B'=B$ which shows why it cannot be a complement for
$B^{\ro(G)}\neq 0$.
\end{exam}

\section{Proof of Theorem~\ref{T:T-FB-general} }  \label{S:DG}
\subsection{Guichardet\texorpdfstring{: $(F_{B})\Longrightarrow(T_{B})$}{}}
\begin{proof}
Assume $G$ does not have $(T_E)$, where $E$ is a Banach space,
and let $\ro:G\to\Orth(E)$ be a representation such that $E/E^{\ro(G)}$
admits almost invariant vectors.
In order to show that $H^{1}(G,\ro)\neq\{0\}$ it suffices to prove that
$B^{1}(G,\ro)\se Z^{1}(G,\ro)$ is not closed.

As was mentioned in Section~\ref{S:prelims}
the space of $\ro$\hy cocycles $Z^{1}(G,\ro)$ is always a Fr\'echet
space (and even a Banach space if $G$ is compactly generated).
Note that $B^{1}(G,\ro)$ is the image of the bounded linear map
\[
    \tau :E\to Z^{1}(G,\ro),\qquad \left(\tau(v)\right) (g)=v-\ro^{\prime}(g)v.
\]
If $\tau(E)$ were closed, and hence a Fr\'echet space,
the open mapping theorem
would imply that $\tau^{-1}:B^{1}(G,\ro)\to E/E^{\ro(G)}$ is a bounded map.
That would mean that for some $M<\infty$ and a compact $K\se G$
\[
    \|v\|\le M\cdot \|\tau(v)\|_{K}=M\cdot \sup_{g\in K}\|\ro(g)v-v\|,
    \qquad v\in E/E^{\ro(G)}
\]
contrary to the assumption that $\ro$ almost contains invariant
vectors.
\end{proof}

\subsection{\texorpdfstring{$(T)\Longrightarrow (F_{L^{p}})$, $0<p\le 2$}{}}

\begin{proof}
Let $G$ be a locally compact group with Kazhdan's property~$(T)$
acting by affine isometries on a closed subspace $B\subseteq
L^{p}(\mu)$ with $0<p \leq 2$. Using Proposition~\ref{F:pos-def} and a slight
modification of a Delorme--Guichardet argument for $(T)\Rightarrow (FH)$
we shall prove that such an action has bounded orbits.
For $1<p\le 2$ uniform convexity of $B\subseteq L^{p}(\mu)$
yields a $G$ fixed point using Lemma~\ref{L:bounded}.

\medskip

Proposition~\ref{F:pos-def} allows to define a family, indexed by
$s>0$, of Hilbert space $\cH_{s}$, embeddings
$U_{s}:B\to\sphere (\cH_{s})$ and unitary representations
$\pi_{s}:G\to\Orth(\cH_{s})$ with the following
properties:
\begin{enumerate}
\item
    The image $U_{s}(B)$ spans a dense subspace of $\cH_{s}$;
\item
    $\ip{U_{s}(x)}{U_{s}(y)}=e^{-s\cdot \|x-y\|^{p}}$ for all $x,y \in B$;
\item
    $U_{s}(gx)=\pi_{s}(g)U_{s}(x)$ for all $x\in B$, $g\in G$.
\end{enumerate}
Indeed, one constructs $\cH_{s}$ as the
completion of the pre-Hilbert space whose vectors are finite
linear combinations $\sum a_{i} x_{i}$ of points $x_{i}\in B$,
and the inner product is given by
\[
    \ip{\sum a_{i} x_{i}}{\sum b_{j} y_{j}}
    =\sum_{i,j} a_{i}\overline{b_{j}} e^{-s\cdot \|x_{i}-y_{j}\|^{p}}.
\]
The representation $\pi_{s}$ can be constructed (and is uniquely
determined) by property~(3).

Since $G$ is assumed to have Kazhdan's property~$(T)$, for some
compact subset $K\se G$ and $\gep>0$, any unitary $G$\hy
representation with $(K,\gep)$\hy almost invariant vectors has a
non-trivial invariant vector.

Let $x_{0}\in B$ be fixed. The isometric $G$\hy action is
continuous, so $K x_{0}$ is a compact and hence bounded subset of
$B$, hence:
\[
    R_{0}=\sup_{g\in K}\|g x_{0}-x_{0}\|<\infty.
\]
For the unit vectors $u_{s}=U_{s}(x_{0})\in\sphere(\cH_{s})$
we have
\[
    \min_{g\in K}|\ip{\pi_{s}(g)u_{s}}{u_{s}}|\ge e^{-s R_{0}^{p}}\to 1
    \qquad\text{as}\quad s\to 0.
\]
In particular for a sufficiently small $s>0$, $\max_{g\in
K}\|\pi_{s}(g)u_{s}-u_{s}\|<\gep$. Let us fix such an $s$, and
rely on property~$(T)$ to deduce that $\pi_{s}$ has an invariant
vector $v\in\sphere(\cH_{s})$.

We claim that $G$ must have bounded orbits for its
affine isometric action on $B$.
Indeed, otherwise there would exist a sequence $g_{n}\in G$
so that
\[
    \|g_{n}x -y\|\to\infty\qquad\text{and hence}\qquad
    \ip{\pi_{s}(g_{n})U_{s}(x)}{U_{s}(y)}\to 0
\]
for all $x,y\in B$.
This implies that $\langle\pi_s(g_n)w,u\rangle\to 0$ for any
$w,u\in\text{span}(U_s(B))$, and since $\text{span}(U_s(B))$ is dense
in $\cH_{s}$, for any $w,u\in\cH_{s}$. Taking
$w=u=v$, we get a contradiction. Therefore the affine isometric $G$\hy
action on $B$ has bounded orbits, and hence fixes a point in case of
$1<p\le 2$.
\end{proof}

\subsection{Fisher--Margulis\texorpdfstring{: $(T)\Longrightarrow (F_{L^{p}})$, $p<2+\e(G)$}{}}\label{S:FM}

Let $G$ have Kazhdan's property~$(T)$. Fix a compact generating
subset $K$ of $G$.
\begin{lem}
There exists a constant $C<\infty$ and $\e>0$
such that for any $G$-action by affine isometries on a closed subspace
$B\subseteq L^{p}(\mu)$ with $p\in (2-\e,2+\e)$
and any $x\in B$ there exists a point $y\in B$ with
\[
    \|x-y\|\le C\cdot \diam(K\cdot x),
    \qquad
        \diam(K\cdot y)< \frac{\diam(K\cdot x)}{2}.
\]
\end{lem}
\begin{proof}
By contradiction there exists a sequence of subspaces
$B_{n}\subseteq L^{p_{n}}$ with $p_{n}\to 2$,
affine isometric $G$\hy actions on $B_{n}$ and points $x_{n}\in B_{n}$
so that, after a rescaling to achieve $\diam(K\cdot x_{n})=1$, we have
\begin{equation} \label{e:displacement}
    \diam (K\cdot y)\ge \frac{1}{2}\qquad \forall\, y\in\ball(x_{n},n).
\end{equation}
Passing to an ultraproduct of the spaces $B_{n}$ with the marked points $x_{n}$
and the corresponding $G$-actions, one obtains an isometric (hence also affine) $G$-action
on a Hilbert space $\cH$, because the limit of
$L^{p}$\hy parallelogram as $p\to 2$ is the parallelogram identity,
which characterizes Hilbert spaces. (The action is well-defined because $K$
generates $G$ and we ensured $\diam(K\cdot x_{n})=1$.)
If $G$ is a topological group, one needs to ensure continuity of the limit
action by selecting uniformly $K$-equicontinuous sets of vectors
(as in~\cite[6.3]{Shalom}; compare also~\cite{Cherix-Cowling-Straub}).
Due to~\eqref{e:displacement} this $G$\hy action has no fixed points,
contradicting property~$(FH)$ and hence $(T)$ of $G$.
\end{proof}

\begin{proof}[Proof of $(F_{B})$ for $B\subseteq L^{p}(\mu)$, $2\le p<2+\e(G)$]
Now consider an arbitrary affine isometric $G$-action on a closed
subspace $B\subseteq L^{p}$ with $|p-2|<\e$
where $\e=\e(G)>0$ is as in the lemma.
Define a sequence $x_{n}\in B$ inductively, starting from an arbitrary $x_{0}$.
Given $x_{n}$, let $R_{n}=\diam(K\cdot x_{n})$. Then applying the lemma
there exists $x_{n+1}$ within the ball $\ball(x_{n},C\cdot R_{n})$ so that
\[
    \diam(K\cdot x_{n+1})< R_{n}/2.
\]
We get $R_{n}<R_{0}/2^{n}$ and $\sum \|x_{n+1}-x_{n}\|<\infty$.
The limit of the Cauchy sequence $\{x_{n}\}$ is a $G$-fixed point.
\end{proof}

\medskip

\begin{ques}
For a given group $G$ with property~$(T)$, what can be said about the following invariant?
$$p(G) :=\inf\big\{ p: \text{ $G$ fails to have~$(F_B)$ for some closed subspace
$B\subseteq L^{p}$ }\big\}.$$
For instance, Pansu's aforementioned result~\cite{Pansu} shows that $p(G)\leq 4n+2$ for $G=\Sp_{n,1}(\BR)$.
\end{ques}

\section{Proof of Theorem~\ref{T:TLp}}  \label{S:TandTB}

We start with the first assertion of the theorem: $(T)\Rightarrow
(T_{B})$ for $B$ being an $L^{p}$\hy related space as in (i), (ii)
or (iii) in the theorem. We first reduce to the case (i) where
$B=L^{p}(\mu)$ with $1\leq p<\infty$. Then using
Corollary~\ref{C:Hardin} of Hardin's extension theorem,
$(T_{L^{p}(\mu)})$ implies~$(T_{B})$ for subspaces $B\se
L^{p}(\mu)$ where $p\neq 4,6,\dots$ as in (ii), and the duality
argument (Corollary~\ref{C:TB<=>TB}) gives the result for
quotients of $L^{q}(\mu)$ with $q\neq 4/3, 6/5, \dots$ as in
(iii). Hence it suffices to prove~$(T)\Rightarrow
(T_{L^{p}(\mu)})$ for $1\leq p<\infty$. We give two proofs for
this implication.

Let us note that our restriction on $p$ and $q$ when taking
subspaces/quotients comes from our use of Hardin's theorem.

\begin{ques}
Does property~$(T)$ implies property~$(T_B)$ for any closed
subspace and any quotient $B$ of $L^p(\mu )$ for any $1<p<\infty$~?
\end{ques}

\subsection{Property~\texorpdfstring{$(T)$ Implies~$(T_{L^{p}(\mu)})$}{(T) Implies (TLp)} -- First Proof}
\label{SS:T=>TLp1}
\begin{proof}
Assuming that a locally compact group $G$ fails to have property
$(T_{L^{p}(\mu)})$ for some $1\leq p<\infty$, we are going to show
that $G$ does not have (T). We may and will assume $p\neq 2$; write
$B=L^p(\mu)$ and $\cH=L^2(\mu)$.
Using Remark~\ref{R:TB} there is a representation
$\ro:G\to\Orth(B)$ so that for the canonical complement
$B^{\prime}$ of $B^{\ro(G)}$ the restriction
$\ro^{\prime}:G\to\Orth(B^{\prime})$ almost has invariant vectors,
i.e. there exist unit vectors  $v_{n}\in \sphere(B^{\prime})$ so
that
\[
    f_{n}(g)=\|\ro(g)v_{n}-v_{n}\|
\]
converges to $0$ uniformly on compact subsets of $G$.

We shall obtain a related unitary, or orthogonal, representation
$\pi:G\to\Orth(\cH)$ using the following:
\begin{lem}\label{L:conj-by-Mazur}
For $p\neq 2$, the conjugation $U\mapsto M_{p,2}\circ U \circ
M_{2,p}$ by the non-linear Mazur map sends $\Orth(B)$ to
$\Orth(\cH)$.
\end{lem}
\begin{proof}
Follows from Banach--Lamperti description of $\Orth(B)$
(Theorem~\ref{T:Banach-Lamperti}) by calculation.
\end{proof}
Let us then define $\pi:G\to\Orth(\cH)$ by
$\pi(g)=M_{p,2}\circ \ro(g) \circ M_{2,p}$.
Note that $M_{p,2}$ maps $B^{\ro(G)}$ onto $\cH^{\pi(G)}$.

As $\sphere(B^{\prime})$ is uniformly separated (in fact is at distance $1$)
from $B^{\ro(G)}$, the uniform continuity of the Mazur map (Theorem~\ref{T:Mazur})
implies that $u_{n}=M_{p,2}(v_{n})$ is a sequence in $\sphere(\cH)$
such that  $\dist(u_{n},\cH^{\pi(G)})\ge\delta>0$ and
$\varphi_{n}(g)=\|\pi(g)u_{n}-u_{n}\|\to 0$
uniformly on compact subsets of $G$.
Let $w_{n}$ denote the projections of $u_{n}$ to
$\cH=(\cH^{\pi(G)})^{\perp}$.
Then
\[
    \|w_{n}\|\ge\delta>0\qquad\text{and}\qquad\|\pi(g)w_{n}-w_{n}\|\le \varphi_{n}(g)\to 0
\]
uniformly on compacta.
Thus the restriction $\pi^{\prime}$ of $\pi$ to $\cH^{\prime}$
does not have $G$\hy invariant vectors, but almost does.
Hence $G$ does not have Kazhdan's property~$(T)$.
\end{proof}

\begin{rem}\label{rem:T_Lp:specific}
In fact, the above proof has established the following more specific statement. Let $G$ act measurably on a $\sigma$-finite measure space. Denote by $\ro_p$ the associated linear isometric representation on $L^p$, namely the quasi-regular representation twisted by the $p$-th root of the Radon--Nikod\'ym derivative. Then, the existence of almost invariant vectors in $L^p \big/ (L^p)^{\ro_p(G)}$ is independent of $1\leq p <\infty$.
\end{rem}

\subsection{Property~\texorpdfstring{$(T)$ Implies~$(T_{L^{p}(\mu)})$}{(T) Implies (TLp)} -- Second Proof}
\label{SS:T=>TLp2}

\begin{proof}
For $1<p\le 2$ we have $(T)\Rightarrow (F_{L^{p}(\mu)})\Rightarrow (T_{L^{p}(\mu)})$
by Theorem~\ref{T:T-FB-general} (1) and (2).
Using duality (Corollary~\ref{C:TB<=>TB}) this implication extends
to $L^{p}(\mu)$ with $2<p<\infty$.
\end{proof}

\subsection{Property~\texorpdfstring{$(T_{L^{p}})$ Implies~$(T)$}{(TLp) Implies (T)}}\label{SS:TLp=>T}

\begin{proof}
Assume that $G$ is not Kazhdan, i.e. $G$ admits a unitary representation
$\pi$ almost containing (but not actually containing) non-trivial
invariant vectors. Connes and Weiss~\cite{CW} showed how to find such
a representation of the form $L^{2}_{0}(\mu)$. More precisely,
 they construct a measure-preserving, \emph{ergodic} $G$\hy action on a probability space $(X,\mu)$
which admits a a sequence $\{E_{n}\}$ of \emph{asymptotically
invariant} measurable subsets, namely
\begin{equation}\label{e:ai}
   \forall\, g\in G \qquad \mu(g E_{n}\bigtriangleup E_{n})\to 0
   \qquad\text{whilst}\qquad
   \mu (E_n)=1/2.
\end{equation}
Consider the unitary $G$\hy representation $\pi^{\prime}$ on
$L^{2}_{0}(\mu)$ -- the space of \emph{zero mean} square
integrable functions, which is the orthogonal complement of the
constants. Then $\pi^{\prime}$ does not have non-trivial invariant
vectors because of ergodicity; but it almost does, namely $f_n
=2\cdot{\bf 1}_{E_{n}}-1$.

For a given $1\leq p<\infty$, consider the linear isometric $G$\hy
representation $\ro$ on $B=L^{p}(\mu)$, $\ro(g)f(x)=f(g^{-1}x)$.
Then $B^{\ro(G)}=\BR\,{\bf 1}$ -- the constants, and its canonical
complement is
\[
    B^{\prime}=L^{p}_{0}(\mu)=\Set{f\in L^{p}(\mu)}{\int f\,d\mu=0}.
\]
The above sequence $\{f_{n}\}$ lies in $L^{p}_{0}(\mu)$,
consists of unit vectors and still
satisfies $\|\ro(g)f_{n}-f_{n}\|_p\to 0$.
Hence failing to have Kazhdan's
property~$(T)$ a group $G$ does not have~$(T_{L^{p}(\mu)})$
either.

In the original paper~\cite{CW}, Connes and Weiss considered
discrete groups.
In a similar context the case of locally compact groups was also
considered by Glasner and Weiss (see~\cite[Section~3]{GW}
and references therein).
One way to treat the non-discrete case, is the following:
start from a unitary representation $\pi$ of a given lcsc $G$
which has almost invariant vectors but no invariant ones, and
apply the original Connes--Weiss Gaussian construction to
the restriction $\pi|_{\Gamma}$ of $\pi$ to some dense countable
subgroup $\Gamma\se G$.
This gives an ergodic measure-preserving $\Gamma$\hy action
on a probability space $(X,\mu)$ with an asymptotically invariant
sequence $\{E_{n}\}$ on $X$.
The fact that the representation $\pi|_{\Gamma}$ came
from $G$ is manifested by the fact that it is continuous in the
topology on $\Gamma$ induced from $G$.
It can be shown to imply that the $\Gamma$\hy representation on
$L^{2}_{0}(X,\mu)$ is also continuous, hence extends to $G$, and thus
the $\Gamma$\hy action on $(X,\mu)$ extends to a measurable $G$\hy action.
This construction gives a uniform convergence in~\eqref{e:ai} on compact subsets
of $G$.
\end{proof}

\section{Fixed Point Property for Higher Rank Groups}\label{S:higher-rank}

\subsection{}
The objective of this section is to prove Theorem~\ref{T:higherrankFLp};
we start with some preliminaries for the \emph{linear} part.
\medskip

The first ingredient needed for the proof is an analogue of
Howe--Moore's theorem on vanishing of matrix coefficients, or
rather its corollary analogous to Moore's ergodicity theorem,
extended to the framework of uniformly equicontinuous representations
on superreflexive Banach spaces. The ucus Banach space version of Howe--Moore
is due to Yehuda Shalom (unpublished). With his kind permission we have included the
argument in Appendix~\ref{Appendix}. Here we shall use the following
corollary, which we formulate for the case of simple groups.
\begin{cor}[Banach space analogue of Moore's theorem]
\label{T:Moore} Let $k$ be a local field and let $G={\mathbf G}(k)$
be the $k$\hy points of a Zariski connected isotropic simple
$k$\hy algebraic group $\mathbf G$.
Let $G^+$ be the image of the simply connected form $\tilde{G}$ in $G$ under the cover map.
Let $H\se G^+$ be a closed
non-compact subgroup.

\nobreak
Then for any superreflexive space $B$ and any
continuous uniformly equicontinuous linear $G$\hy representation $\ro:G^+\to\GL(B)$,
$B^{\ro(H)}=B^{\ro(G^+)}$ and the canonical complements
with respect to both $\ro(G^+)$ and $\ro(H)$ coincide,
and can be denoted just by $B^{\prime}$.
\end{cor}

\begin{proof}
By Proposition~\ref{P:superreflex-ucus}, we may assume that $B$ is a ucus Banach space
and $\ro$ is a linear isometric representation $\ro:G\to\Orth(B)$. Now the statement follows
readily from Theorem~\ref{T:HM}.
\end{proof}

\subsection{}
The second ingredient is \emph{strong relative property~$(T)$}.
It will be used to prove Claim~\ref{Claim:AxH}
below which is the only part which is specific to $L^{p}$\hy like spaces.
The rest of the argument applies to all affine isometric actions on ucus Banach spaces,
or all uniformly equicontinuous affine actions on a superreflexive space.
\begin{defn}\label{D:relT}
Let $H\ltimes U$ be a semi-direct product of locally compact groups.
We shall say that it has
\begin{description}
\item[strong relative property~$(T)$]
    if for any unitary representation $\pi$ of $H\ltimes U$ for which
    $H$ almost has non-trivial invariant vectors, $U$ has invariant vectors.
\item[strong relative property~$(T_{B})$]
    where $B$ is a Banach space, if for any linear isometric
    representation $\ro:H\ltimes U\to\Orth(B)$ the linear isometric
    $H$\hy representation $\ro^{\prime}:H\to\Orth(B/B^{\ro(U)})$
    does not almost have non-trivial invariant vectors.
\end{description}
\end{defn}
\begin{rems}\label{R:relT}~
\begin{enumerate}
\item
    The first definition is a variant of ``relative property~$(T)$''.
    The latter usually refers to a pair
    of groups $G_{0}\se G$ and requires that any unitary $G$\hy representation
    with $G$\hy almost invariant vectors, has non-trivial $G_{0}$\hy invariant vectors.
    Strong relative property~$(T)$ for $H\ltimes U$ implies, but is not equivalent to,
    relative property~$(T)$ for $(H\ltimes U, U)$.
    In fact $\SL_{2}(\BR)\ltimes \BR^{2}$ has the strong relative~$(T)$ and thus relative (T)
    as well, whilst its lattice $\SL_{2}(\BZ)\ltimes \BZ^{2}$
    does not have strong relative~$(T)$ even though the pair
    $(\SL_{2}(\BZ)\ltimes \BZ^{2},\BZ^{2})$ has relative property~$(T)$.
    (For the latter, cf. M. Burger's appendix in ~\cite{HV}. For the former, consider
    the representation on $\ell^2(\BZ^2)$ induced by the affine action on $\BZ^2$.)
\item
    If $B$ is a ucus Banach space, then the canonical splitting with respect
    to $\ro(U)$, namely $B=B^{\ro(U)}\oplus B^{\prime}$ is preserved
    by $\ro(H)$ which normalizes $\rho (U)$ (Corollary~\ref{normelizers}).
    Hence, as in Remark~\ref{R:TB}, for ucus space $B$
    strong relative property~$(T_{B})$
    requires that the \emph{restriction} of $\ro(H)$ to $B^{\prime}$ does
    not almost have invariant vectors.
    Strong relative~$(T_{\cH})$ for a Hilbert space $\cH$
    is equivalent to the strong relative~$(T)$.
\end{enumerate}
\end{rems}
\begin{lem} \label{L:strong-relative-T}
A semi-direct product $H\ltimes U$ with strong relative
property~$(T)$ has strong relative property~$(T_{B})$ for all
$L^{p}$\hy related Banach spaces $B$ of types {\rm (i), (ii),
(iii)} as in Theorem~\ref{T:TLp}.
\end{lem}
\begin{proof}
This is analogous to the proof of $(T)\Rightarrow (T_{B})$ given
in Section~\ref{SS:T=>TLp1}. First observe that the extension
Theorem~\ref{C:Hardin} and a duality argument (based on
Proposition~\ref{stam}) reduce the statement to the
case (i) of $B=L^{p}(\mu)$.

Thus we assume that $B=L^{p}(\mu)$ with $p\neq 2$, and
$\ro:H\ltimes U\to\Orth(B)$ is a linear isometric representation.
Let $B=B^{\ro(U)}\oplus B^{\prime}$ be the canonical splitting
with respect to $U$. It is preserved by $\ro(H)$ because $H$ normalizes $U$.
Now let $\pi=M_{p,2}\circ \ro\circ M_{2,p}$ be the conjugate of $\ro$ by the Mazur map.
Then $\pi$ is an orthogonal representation $\pi:H\ltimes U\to \Orth(\cH)$
where $\cH=L^{2}(\mu)$ (Lemma~\ref{L:conj-by-Mazur}).

If $H\ltimes U$ fails to have strong relative~$(T_{B})$, then there exist
$x_{n}\in \sphere(B^{\prime})$ so that $\|\ro(h)x_{n}-x_{n}\|\to 0$
uniformly on compact subsets of $H$.
Uniform continuity of $M_{p,2}$ and the fact that
$\dist(\sphere(B^{\prime}),\sphere(B^{\ro(U)}))=1$, imply
that for $v_{n}=M_{p,2}(x_{n})$
\[
    \dist(v_{n},\cH^{\ro(U)})\ge\delta>0\qquad \|\pi(h)v_{n}-v_{n}\|\to0
\]
uniformly on compact subsets of $H$.
Taking projections of $v_{n}$ to $\cH^{\prime}$ we show that
in this case $H\ltimes U$ does not have strong relative property~$(T)$.
\end{proof}
%
\subsection{Proof of Theorem~\ref{T:higherrankFLp}}
We first show that we can assume that $G$ is connected and simply connected.
Assuming that Theorem~\ref{T:higherrankFLp} is known for
$\tilde{G_0}$ and lattices therein; we will prove it for $G$ and its lattices. 
For any affine isometric
action of $G$ on $B$ there is an associated action of $\tilde{G_0}$,
inflated via the covering map $\tilde{G_0}\rightarrow G$.
$\tilde{G_0}$ has a fixed point by assumption, hence $G$ has a
compact orbit, as the cokernel of the covering map is
compact~\cite[Theorem~2.3.1(b)]{Margulis}. 
It follows that $G$ has a fixed point as well.
A similar argument applies to lattices:
For a given lattice $\Gamma$ in $G$ its $\tilde{\Gamma}$ by the covering 
map is a lattice in $\tilde{G_0}$, and its projection is of finite index in $\Gamma$.
Every affine isometric action of $\Gamma$gives rise to an affine isometric 
action of $\tilde{\Gamma}$,
which, by assumption, has a fixed point.
It follows that $\Gamma$ has a finite orbit, and therefore fixes a point.

Hereafter we will assume that $G$ is connected and simply connected.
In that case $G$ decomposes into a direct product of simply connected
almost simple groups
$G=\prod G_{i}$ \cite[Proposition 1.4.10]{Margulis}.

In view of (the independent) Sections~\ref{sec_ind_gen} and~\ref{sec_induction}, more specifically
Proposition~\ref{P:lattices}(2) and the discussion following Definition~\ref{defn_p_int},
property~$(F_{B})$ for $G=\prod G_{i}$ is inherited by its lattices.
Thus it suffices to consider the ambient group $G=\prod G_{i}$ only.
By Proposition~\ref{P:generalities}(3)
the statement reduces to that about almost-simple factors $G_{i}$.

So we are left proving the theorem for $G={\mathbf G}(k)$, a higher rank connected, simply-connected, almost-simple group.
Using Proposition~\ref{P:superreflex-affine}, we assume that $B$ is a ucus Banach space
and we consider a $G$\hy action on $B$ by affine isometries, with
$\ro:G\to\Orth(B)$ denoting the linear part of the action.
Let $B=B^{\ro(G)}\oplus B^{\prime}$ be the canonical
decomposition and $\ro^{\prime}:G\to\Orth(B^{\prime})$ denote
the corresponding sub-representation.

\begin{clm}[For $L^{p}$\hy like spaces]\label{Claim:AxH}
$G$ contains a direct product $A\times H$ so that
\begin{enumerate}
\item
    The restriction $\ro^{\prime}|_{H}:H\to\Orth(B^{\prime})$ does
    not almost contain invariant vectors.
\item
    $A$ contains a
    non-trivial semisimple element, and in particular it is not
    compact.
\end{enumerate}
\end{clm}
\begin{proof}
Any higher rank almost-simple group $G={\mathbf G}(k)$ is known to contain
a subgroup whose simply-connected cover is isomorphic to either $G_{0}=\SL_{3}(k)$ or $G_{0}=\Sp_{4}(k)$
\cite[Theorem~1.6.2]{Margulis}.
In the first case $G_{0}=\SL_{3}(k)$ contains the semi-direct product
$H_{0}\ltimes U_{0}=\SL_{2}(k)\ltimes k^{2}$ embedded in $\SL_{3}(k)$ as
\[
    \left\{\left.
    \left(\begin{array}{ccc}
    a & b& x\\
    c & d & y\\
    0 & 0 & 1\end{array} \right) \right|
    ad-bc=1 \right\}
\]
where $U_{0}\cong k^{2}$ is the subgroup given by $a=d=1$,
$b=c=0$. It is normalized by the copy $H_{0}$ of $\SL_{2}(k)$
embedded in the upper left corner. Let $A_{0}\se\SL_{3}(k)$ be
the subgroup $\diag[\lambda,\lambda,\lambda^{-2}]$, $\lambda\in
k^{*}$, which centralizes $H_{0}$ in $G_{0}$, and let $A$ and
$H\ltimes U$ denote the corresponding subgroups in $G$.

The semi-direct product $\SL_{2}(k)\ltimes k^{2}$ is known to have
strong relative property~$(T)$. Hence it has strong relative
property~$(T_{B})$ for $L^{p}$\hy related spaces $B$
(Lemma~\ref{L:strong-relative-T}). By~\ref{T:Moore} we have
$B^{\ro(G)}=B^{\ro(U)}$ and we have denoted by $B^{\prime}$  the
common canonical complement. Then (1) follows from the strong
relative property~$(T_B)$ for $H\ltimes U$, while (2) is clear
from the construction.

In the second case $G$ contains a copy of $G_{0}=\Sp_{4}(k)$ which
is usually defined as a subgroup of $\SL_{4}(k)$ by
\[
    \Sp_{4}(k)=\left\{ g\in \SL_{4}(k) \,\mid\, ^{t}g J g= J\right\},\qquad \text{where}\qquad J=\left(\begin{array}{cc} 0 & I \\ -I & 0\end{array}\right).
\]
The semi-direct product $H_{0}\ltimes U_{0}$ embedded in $\SL_{4}(k)$ is
\[
    \left\{\left.  \left(\begin{array}{cc}
        A & B\\ 0 & ^{t}A^{-1}
    \end{array} \right) \right|
    A\in\SL_{2}(k),\ ^{t}B=A^{-1} B\, (^{t}A) \right\}
\]
with $H_{0}$ denoting the image $A\mapsto \diag[A, ^{t}A^{-1}]$ of $\SL_{2}(k)$,
and $U_{0}$ the normal Abelian subgroup
\[
    \left\{\left.  \left(\begin{array}{cc}  I & B\\ 0 & I\end{array} \right)\,
    \right|\,  ^{t}B=B \right\}.
\]
The semi-direct product $H_{0}\ltimes U_{0}$ actually lies in
$\Sp_{4}(k)$, it is isomorphic to $\SL_{2}(k)\ltimes S^{2}(k)$,
where $S^{2}(k)$ is the space of symmetric bilinear forms on
$k^{2}$ with the natural $\SL_{2}(k)$ action. This semi-direct
product is also known to have strong relative property~$(T)$, and
therefore strong relative~$(T_{B})$. $H_{0}$ is centralized by
$A_{0}=\Set{\diag[\lambda,\lambda,\lambda^{-1},\lambda^{-1}]}{\lambda\in
k^{*}}$. As in the $G_{0}=\SL_{3}(k)$ case, we conclude that the
corresponding product $A\times H\se G$ satisfies (1) and (2).
The claim is proved.
\end{proof}

We now turn to the affine isometric $G$\hy action defined by a
$\ro$\hy cocycle $c\in Z^{1}(\ro)$. We shall prove that $c\in
B^{1}(\ro)$ i.e. that $G$ has a global fixed point. Write
$c(g)=c_{0}(g)+c^{\prime}(g)$ with $c_{0}(g)\in B^{\ro(G)}$ and
$c^{\prime}(g)\in B^{\prime}$ where $B=B^{\ro(G)}\oplus
B^{\prime}$ is the canonical splitting. Then $c_{0}:G\to B$ is a
homomorphism into the (additive) Abelian group. As $G$ has compact
Abelianization, $c_{0}(g)\equiv 0$, which means that the affine
$G$\hy action preserves each affine subspace $p+B^{\prime}$. Hence
both the affine $G$\hy action and the representation can be
restricted to $B^{\prime}$.

Claim~\ref{Claim:AxH} provides an input for the following general
lemma:
\begin{lem} \label{L:AxHfix}
Let a direct product of topological groups $A\times H$ act by
affine isometries on a Banach space $B$. Suppose that the
associated linear isometric representation $\ro$ restricted to $H$
does not almost have invariant vectors. Then the affine action of
$A$ has bounded orbits in $B$. In particular, if $B$ is uniformly
convex, then $A$ has a fixed point in $B$.
\end{lem}
\begin{rem}
In the uniformly convex case, this follows of course from the stronger splitting theorem
(Theorem~\ref{thm_splitting}); compare also with Theorem~\ref{thm_product}
below for the weaker assumption that
the \emph{product} does not almost have invariant vectors.
\end{rem}
\begin{proof}[Proof of the lemma]
Let $\ro:A\times H\to\Orth(B)$ and $c\in Z^{1}(\ro)$ denote the
associated linear isometric representation and the
translation cocycle.
The commutation relation between any $h\in H$ and $a\in A$ gives
\[
    c(h)+\ro (h) c(a) =c(ha)=c(ah)=c(a)+\ro (a) c(h)
\]
which can be rewritten as
\[
    (I-\ro(h)) c(a)=(I-\ro(a)) c(h).
\]
By the assumption on $\ro(H)$, there exists a compact subset $K\se H$
and an $\gep>0$ so that
$\max_{h\in K} \|\ro(h)v-v\|\ge \gep\cdot \|v\|$ for all $v\in B$.
Let $R=\max_{h\in K} \|c(h)\|<\infty$.
Then for $a\in A$
\[
    \gep\cdot \|c(a)\|\le \max_{h\in K}\left\|(I-\ro(h))c(a)\right\|
    \le 2 R.
\]
Hence $\sup_{a\in A}\|c(a)\|\le 2R/\gep$, i.e. the $A$\hy orbit of $0$ is bounded.
If $B$ is uniformly convex then the circumcentre of this orbit is an $A$\hy fixed point
as in Lemma~\ref{L:bounded}.
\end{proof}

Let us restrict the $G$\hy action to $B^{\prime}$. For $g\in G$
let $\Fix(g)$ denote the set of $g$\hy fixed points in
$B^{\prime}$. It follows from Claim~\ref{Claim:AxH} and
Lemma~\ref{L:AxHfix} that for some non-elliptic semisimple element
$a\in A\se G$, $\Fix(a)\subseteq B^{\prime}$ is non-empty.

\smallskip

At this point we propose two different ways to conclude the proof.

\smallskip

\noindent
\textbf{First argument.~}For any $g$ we define
$$U(g) = \big\{h\in G : \lim_{n\to\infty} g^{-n} h g^n = e\big\}.$$
In analogy to the \emph{Mautner phenomenon}~\cite[II.3.2]{Margulis}, we remark that in any continuous $G$-action by isometries on any metric space, every $g$-fixed point is $U(g)$-fixed. Indeed, if $gx=x$ then $d(hx, x) = d( g^{-n} h g^n x, x) \to 0$. It follows that any $g$-fixed point is fixed by both $U(g)$ and $U(g^{-1})$.  Back to the setting of Theorem~\ref{T:higherrankFLp}, we apply this to $a$ and find that any $a$-fixed point is fixed by $G$ because the latter is generated by $U(a)\cup U(a^{-1})$, see~\cite[I.1.5.4(iii)]{Margulis}.

\smallskip

\noindent
\textbf{Second argument.~}Note that $\Fix(a)$ cannot contain more than one point. Indeed, if
$x,y\in\Fix(a)$ then $x-y\in B^{\prime}$ is a fixed vector for the
linear isometry $\ro(a)$. Since the cyclic group $\langle a
\rangle$ is unbounded, the ucus analogue of Moore's
ergodicity~\ref{T:Moore} implies that $B^{\ro(\langle
g\rangle)}=B^{\ro(G)}$. Hence $x-y=0$. Thus $\Fix (a)=\{ x_0\}$.

If $g,h\in G$ commute then $\Fix(g)$ is an $h$\hy invariant set.
So if $\langle g\rangle $ is unbounded and $\Fix(g)\neq\emptyset$
then $\Fix(g)$ is a single point fixed by $h$. If $\langle
h\rangle$ is also unbounded then $\Fix(g)=\Fix(h)$.
Hence the following lemma implies that $G\cdot x_0=x_0$ and again finishes the proof of Theorem~\ref{T:higherrankFLp}.

\begin{lem}\label{lem:5.8}
For any two non-elliptic semisimple elements $g,h\in G$ there is a
chain $g=g_1,g_2,\ldots,g_n=h$ of non-elliptic semisimple
elements, each commuting with its successor. Furthermore, the set
of non-elliptic semisimple elements generates $G$.
\end{lem}
\begin{proof}
The first claim is equivalent to the connectivity of the Tits
boundary, which is equivalent to the assumption on the rank. Both claims are well-known.
\end{proof}

\begin{rem}
Observe that the somewhat restrictive assumption that $B$ is an $L^{p}$\hy related
space is used only in the proof of Claim~\ref{Claim:AxH}, the rest
of the argument being in the context of general ucus Banach spaces.
\end{rem}
%


\section{Minimal Sets}\label{sec_min}

Let $B$ be a strictly convex reflexive Banach space and $G$ a
group acting on $B$ by affine isometries. Consider the ordered category
$\mC$ of non-empty closed convex $G$\hy invariant subsets of $B$
endowed with $G$\hy equivariant isometric maps and inclusion
order. The goal of this section is to study minimal elements of
$\mC$ (regardless of whether they exist). In
Section~\ref{sec_product} we shall prove their existence, under
certain conditions (see Corollary~\ref{cor_min_exists}).

The Mazur--Ulam theorem states that a surjective isometry between
(real) Banach spaces is affine. It is not known (and probably not true
under no further assumptions) whether the analogous of the
Mazur--Ulam theorem holds in the general context of convex subsets
of Banach spaces. However, for subsets of strictly convex spaces
it is obviously true:

\begin{lem}\label{iso-affine}%
Let $C\se B$ be a convex subset. Then every isometric map $C\to B$ is affine.
\end{lem}

\begin{proof}
It is enough to show that for all $x,y\in C$ and every $0< t< 1$ the point $p=t x + (1-t)y$
is determined metrically.
This is true since by strict convexity
$$\ol\ball(x,(1-t)\|x-y\|)\, \cap\, \ol\ball(y,t\|x-y\|) = \{p\}.$$
\end{proof}

In particular the morphisms of $\mC$ are affine.
Another useful geometric property of closed convex sets in $B$ is the
existence of a nearest point projection.

\begin{lem}\label{exist_pi}%
Let $C$ be a non-empty closed convex subset of $B$.
Then for every $x\in B$ there exist a unique point
$\pi_C(x)\in C$ such that $\|x-\pi_C(x)\|=d(x,C)$.
 \end{lem}

 \begin{proof}
The uniqueness follows from strict convexity. By the Hahn-Banach theorem $C$
is weakly closed since it is closed and convex;
therefore, by reflexivity and the Banach--Alao\u{g}lu theorem we have a nested family
$C\cap\ \ol\ball(x,d)$ of weakly compact sets as $d \searrow d(x,C)$;
its intersection yields existence.
\end{proof}

The map $\pi_C:B\rightarrow C$ is called the \emph{nearest point projection} on $C$.
We remark that it is not continuous in general.
It is continuous for uniformly convex Banach spaces and non-expanding for Hilbert spaces.
Still, the distance between a point and its projection is always a $1$\hy Lipschitz function:

\begin{lem}\label{semi-cont}%
Let $C$ be a non-empty closed convex subset of $B$.
Then the function $x\mapsto \|\pi_C(x) - x\|$ from $B$ to $\BR$ is $1$\hy Lipschitz.
\end{lem}

\begin{proof}
For any $x,y\in B$
$$\|\pi_C (x) - x\| \leq \|\pi_C (y) - x\| \leq \|\pi_C (y) - y\| + \|y - x\|.$$
\end{proof}

\begin{lem}\label{no-conv-map}%
If $C\in \mC$ is a minimal element, then any convex $G$\hy invariant continuous
(or lower semi-continuous) function $\gf:C\to\BR$ is constant.
\end{lem}

\begin{proof}
If $\gf$ were to assume two distinct values $s<t$, then $\gf^{-1}((-\infty,s])$
would be a strictly smaller element of $\mC$.
\end{proof}

\begin{lem}\label{pi-affine}%
Let $C,C'\in\mC$ with $C$ minimal. Then the nearest point projection
$\pi=\pi_{C'}|_C:C\rightarrow C'$ is affine.
\end{lem}

\begin{proof}
For every $x,y\in C$ and $t\in [0,1]$, the definition of $\pi$ implies
\begin{multline}\label{eq-pi-affine}%
    \|\pi\big( tx+(1-t)y\big) -\big( tx+(1-t)y\big)\|  \\
    \leq \|\big( t\pi (x)+(1-t)\pi (y)\big) -\big( tx+(1-t)y\big)\|  \\
    \leq  t\|\pi (x)-x\|+(1-t)\|\pi (y)-y\|.
\end{multline}
It follows that the function $C\rightarrow \BR$, $x\mapsto \|\pi (x)-x\|$ is convex.
Clearly it is $G$\hy invariant, and by Lemma~\ref{semi-cont} it is continuous, hence Lemma~\ref{no-conv-map} implies
that $\|\pi (x)-x\|$ is constant on $C$. This constant must be $d(C,C')$; as both the
right-hand side and the left-hand side in~\eqref{eq-pi-affine} equal $d(C,C')$, it follows that
\[
    \|\big( t\pi (x)+(1-t)\pi (y)\big) -\big( tx+(1-t)y\big)\|
    =  \|\pi\big( tx+(1-t)y\big) -\big( tx+(1-t)y\big)\|.
\]
Therefore, by the uniqueness part of Lemma~\ref{exist_pi}, $t\pi (x)+(1-t)\pi (y)$
must be $\pi\big( tx+(1-t)y\big)$.
\end{proof}

\begin{lem}\label{affine-const}%
If $C\in\mathcal{C}$ is minimal and $T:C\rightarrow B$ is a $G$\hy equivariant affine map,
then there exist a $\ro(G)$\hy invariant vector $b\in B$ such that $T(c)=c+b$ for all $c\in C$.
\end{lem}

\begin{proof}
The map $C\rightarrow \BR$, $x\mapsto \|Tx-x\|$ is $G$\hy invariant, continuous and convex, hence by
Lemma~\ref{no-conv-map} it has a constant value $d\geq 0$.
Since $B$ is strictly convex and $C$ is convex, the affine map $\sigma(x) = Tx - x$
from $C$ to the sphere of radius $d$ in $B$ must be constant.
Its value $b=\gs (C)$ is the desired ($\ro (G)$\hy invariant) translation vector.
\end{proof}

\begin{cor}\label{pi-translates}
The map $\pi_C:C\to C'$ from Lemma~\ref{pi-affine} is in fact a translation.
\end{cor}

\begin{cor}\label{iso-trans}%
If $C,C'\in\mathcal{C}$ are minimal, then they are equivariantly isometric.
Moreover, any equivariant isometry $C\to C'$ is a translation by a $\ro(G)$\hy invariant vector.
\end{cor}

\begin{proof}
By Corollary~\ref{pi-translates}, $\pi_{C'}|_C: C\to C'$ is an isometry;
it is $G$\hy equivariant and hence onto by minimality of $C'$.
The second claim follows from Lemma~\ref{iso-affine} and Lemma~\ref{affine-const}.
\end{proof}

\section{Actions of Product Groups and Splitting}\label{sec_product}

\subsection{}
The main goal of this section is to prove Theorem~\ref{thm_splitting}.
By Proposition~\ref{P:superreflex-affine} we may assume the affine
action to be isometric with respect to a ucus norm on a Banach space $B$.
The main step is the following theorem.
\begin{thm}\label{thm_product}%
Let $G=G_1\times G_2$ be a product of topological groups with a continuous
action by affine isometries on a uniformly convex Banach space $B$ without $G$\hy fixed point.
Assume that the associated linear $G$\hy representation $\ro$ does not almost have
non-zero invariant vectors.
Then there exists a non-zero $\ro(G_i)$\hy invariant vector for some $i\in\{1,2\}$.
\end{thm}

The proof of Theorem~\ref{thm_product} uses minimal sets (in analogy to~\cite{Monod}); notice that
we are in the setting of Section~\ref{sec_min}
since uniformly convex spaces are reflexive and strictly convex~\cite[App.~A]{BL}.
More precisely, we show:

\begin{prop}\label{prop_min_exists}%
Let $G$ and $B$ be as above. Then there exists a minimal non-empty
closed convex $G_1$\hy invariant subset in $B$. In fact, any
non-empty closed convex $G_1$\hy invariant subset contains a
minimal such subset.
\end{prop}

\begin{proof}[Proof of Theorem~\ref{thm_product}]
Proposition~\ref{prop_min_exists} provides a minimal non-empty closed convex
$G_1$\hy invariant set $C\se B$.
If there is no non-zero $\ro(G_1)$\hy invariant vector, Lemma~\ref{affine-const}
(applied to $G_1$) shows that $G_2$
fixes every point of $C$. Since $G_1$ preserves $C$ and $G$ has no fixed point,
$C$ cannot consist of a single point.
Picking two distinct points $x,y\in C$ yields the non-zero $\ro(G_2)$\hy invariant vector $x-y$.
\end{proof}

Recall that uniform convexity is characterized by the positivity of the convexity modulus
$\delta$ defined in Section~\ref{SS:Banach}.
Moreover, $\delta$ is a positive, non-decreasing
function which tends to zero at zero. Defining
\[
    \delta^{-1}(t)=\sup\{\gep~:~\delta(\gep) \leq t\},
\]
$\delta^{-1}$ is easily seen to share the same properties. Furthermore, for every
$\gep>0$, $\delta^{-1}\circ\delta(\gep)\geq \gep$.

\begin{proof}[Proof of Proposition~\ref{prop_min_exists}]
Let $C_0\se B$ be any non-empty closed convex $G_1$\hy invariant subset; we will show that $C_0$ contains
a minimal subset (if no initial $C_0$ was prescribed, one may choose $C_0 = B$).

Pick any $p\in C_0$ and let $C_1\se C_0$ be the closed convex hull of the $G_1$\hy orbit of $p$. By Hausdorff's maximal
principle, we can
chose a maximal chain $\mD$ of non-empty closed convex $G_1$\hy invariant subsets of $C_1$. If $b_C:=\pi_C (0)$ is
bounded as $C$ ranges over $\mD$, then for some $R>0$ we have a nested family of non-empty sets
$\ol\ball(0,R)\cap C$ which are weakly compact by reflexivity, Hahn--Banach theorem and Banach--Alao\u{g}lu theorem.
In particular the intersection
$\bigcap\mD$ is non-empty, thus providing a minimal set for $G_1$. Therefore, we may from now on assume for a
contradiction that the (non-decreasing) net $R_C:=\|b_C\|$ is unbounded over $C\in\mD$. Let $\mD'\se \mD$ be the
cofinal segment defined by $R_C> 0$. We will obtain a contradiction by showing that for every compact $K\se G$,
$\diam(\ro(K)\hat{b}_C)$ tends to zero along $C\in\mD'$, where $\hat{b}_C = \frac{b_C}{R_C}$.

Indeed, choose $K_i\se G_i$ compact with $K\se K_1\times K_2$ and
let $L=\max_{g\in K_1\times K_2} \|g\cdot 0\|$. The choice of
$b_C$ implies $g\cdot b_C\neq 0$ and $R_C \leq \|\frac{b_C +
g\cdot b_C}{2}\|$ for all $g\in G$. Therefore, setting $x =
\frac{b_C}{\|g\cdot b_C\|}$, $y = \frac{g\cdot b_C}{\|g b_C\|}$,
the convexity modulus $\delta_{C,g}:= \delta(\|x - y\|)$ gives
\begin{multline*}
R_C \leq \|\frac{1}{2} b_C+\frac{1}{2} g\cdot b_C\| \leq \|\frac{1}{2}(x+y)\|\cdot\|g\cdot b_C\|
 \leq (1-\delta_{C,g})\cdot (\|g\cdot b_C-g\cdot 0\|+\|g\cdot 0\|)\\
\leq (1-\delta_{C,g})(R_C+L) \leq R_C(1+\frac{L}{R_C}-\delta_{C,g})\kern1cm \forall\,g\in K_1.
\end{multline*}
Therefore $\delta_C:=\sup_{g\in K_1}\delta_{C,g}\leq \frac{L}{R_C}\to 0$ along $C\in\mD'$ and hence
$$\sup_{g\in K_1} \frac{\|g\cdot b_C - b_C\|}{\|g\cdot b_C\|}\leq \delta^{-1}(\delta_C) \to 0.$$
Using $\|g\cdot b_C\| \leq \| g\cdot b_C - g 0\| +L \leq R_C + L$, it follows that
\begin{equation}\label{eq_decay_G1}%
\sup_{g\in K_1} \frac{\|g\cdot b_C - b_C\|}{R_C} \to 0 \kern.5cm\text{along}\ C\in\mD'.
\end{equation}
On the other hand, for every $g\in G_2$, the function $z\mapsto
\|g\cdot z - z\|$ is continuous, convex and $G_1$\hy
invariant;therefore, it is bounded by $\|g\cdot p-p\|$ on $C_1$.
Setting $L' = \max_{g\in K_2}\|g\cdot p-p\|$, it follows now that
for all $k=(g_1, g_2)\in K$ we have
\begin{multline*}
R_C \cdot \|\ro(k)\hat{b}_C - \hat{b}_C \| = \|k b_C - b_C - k 0\| \leq \|g_1 b_C - b_C\| + \|g_2 g_1 b_C - g_1 b_C\|+ L\\
\leq \|g_1 b_C - b_C\| + L'+ L.
\end{multline*}
Thus, in view of~\eqref{eq_decay_G1}, $\diam(\ro(K)\hat{b}_C)$ goes to zero as claimed.
\end{proof}

\begin{proof}[Proof of Theorem~\ref{thm_splitting}]
We adopt the notation and assumptions of that theorem; let $\ro$ be the linear part of the action. Assume first $n=2$.
Since we have in particular $B^{\ro(G)}=0$, Corollary~\ref{C:split:ucus} yields a canonical splitting
$B = B^{\ro(G_1)}\oplus B^{\ro(G_2)}\oplus B_0$ invariant under $\ro(G)$. Decomposing the cocycle $G\to B$ along this
splitting shows that up to affine isometry we may assume that the affine $G$\hy space $B$ splits likewise as affine product
of affine spaces with corresponding linear parts. However, Theorem~\ref{thm_product} shows that the resulting affine
$G$\hy action on $B_0$ must have a fixed point since $B_0^{\ro(G_i)}=0$. Therefore we obtain a $G$\hy invariant affine
subspace $G$\hy isometric to $B^{\ro(G_1)}\oplus B^{\ro(G_2)}$ in $B$, as claimed.

In order to obtain the general case $n\geq2$, we only need to observe that Corollary~\ref{C:split:ucus} applied to the
product $G_1\times \prod_{i\geq 2} G_i$ allows us to apply induction on $n$.
\end{proof}

\begin{rem}\label{rem:product:subspace}
The above proof characterizes as follows the subspaces $B_i\se B$
appearing in the statement of Theorem~\ref{thm_splitting}: Upon
possibly replacing the $B_i$ with the corresponding linear
subspace (which corresponds to replacing the cocycles with
cohomologous cocycles), we have $B_i = B^{\ro(G'_i)}$ for $G'_i =
\prod_{j\neq i} G_j$.
\end{rem}

\subsection{A More Geometric Approach to Theorem~\ref{T:higherrankFLp} and a Step Towards Conjecture~\ref{C:higherrankFB}}

Before going on towards the superrigidity theorem, let us explain
a more geometric, and seemingly more general, approach to prove
$(T_B)\Rightarrow (F_B)$, which is based on minimal sets. First we
shall formulate a very general statement in the vein of Conjecture
\ref{C:higherrankFB}:

\begin{thm}\label{thm:TiwF}
Let $B$ be a ucus Banach space and $G$ a topological group with
property~$(T_B)$ and compact Abelianization. Then for any
continuous affine isometric action of $G$ on $B$ there is a
minimal non-empty closed convex subset $C\se B$. Moreover
$\Aut_G(C)$ is trivial, $C\se B'$ and $C$ is unique up to
translations by a $\ro (G)$-invariant vector.
\end{thm}

The proof of Theorem~\ref{thm:TiwF} relies on the following consequence of our discussion of minimal sets:

\begin{cor}\label{cor_min_exists}%
Let $G$ be a topological group with a continuous action by
affine isometries on a uniformly convex Banach space $B$. Assume that the
associated linear representation does not almost have non-zero
invariant vectors. Then there exists a unique minimal non-empty
closed convex $G$\hy invariant subset $C_0\se B$. Moreover, there
are no non-trivial $G$\hy equivariant isometries of $C_0$.
\end{cor}

\begin{rem}\label{rem_min_exists}
In view of the additional statement of
Proposition~\ref{prop_min_exists}, the set $C_0$ is contained in
every non-empty closed convex $G$\hy invariant subset. Thus it is
indeed the (non-empty) intersection of all those subsets.
\end{rem}

\begin{proof}[Proof of Corollary~\ref{cor_min_exists}]
For the existence of $C_0$, we may apply
Proposition~\ref{prop_min_exists} if $G=G_1\times 1$ has no fixed
point, or otherwise take such a fixed point for $C_0$. Both
uniqueness and the additional statement follow now from
Corollary~\ref{iso-trans}.
\end{proof}

\begin{proof}[Proof of Theorem~\ref{thm:TiwF}]
Since $G$ has compact Abelianization, the $\ro(G)$-invariant
subspace $B'$ is in fact $G$-invariant as an affine space,
as the projection of the cocycle to $B^{\ro(G)}$ must be a homomorphism.
It
follows that every minimal non-empty closed convex $G$-invariant
set is contained in some coset of $B'$. The existence and
uniqueness of such subset $C$ inside $B'$ follows from Corollary~\ref{cor_min_exists}.
The fact that any two such sets are
different by a $\ro (G)$-invariant vector is a consequence of
Corollary~\ref{iso-trans}.
\end{proof}

Let us now describe an alternative proof for Theorem~\ref{T:higherrankFLp}.
Let $B$ be an $L^p$-related Banach space as in
Theorem~\ref{T:higherrankFLp}.
We reduce to the case where $G$ is connected, simply-connected and
almost-simple as in Section \ref{S:higher-rank}.
Now $G$ either contains a copy of $\SL_3(k)$ or of $\Sp(4,k)$
that, in each case, contains a semidirect product $H\ltimes U$
with the strong relative property~$(T_B)$ (see
Lemma~\ref{L:strong-relative-T} and the proof of
Claim~\ref{Claim:AxH} ). We decompose $B=B^{\ro (U)}\oplus B'$
according to that $U$ action (note that by Howe-Moore $B^{\ro
(U)}=B^{\ro (G)}$). Then $B'$ is invariant under the affine action
of $G$ and $H$ does not almost has invariant vectors in $B'$.
Hence, by Corollary~\ref{cor_min_exists} there is a unique minimal
non-empty closed convex $H$-invariant subset $C\se B'$ and it has
no non-trivial automorphisms which commute with the $H$-action.
Since, by Claim~\ref{Claim:AxH}, the centralizer of $H$ is
non-compact it follows by Howe--Moore that $|C|=1$.
As in Section \ref{S:higher-rank}, we can now finish the proof
using Lemma \ref{lem:5.8}.

\section{Induction and Superrigidity}\label{S:superrigidity}

Let $\Gamma<G=G_1\times \cdots \times G_n$ be a lattice in a product of $n\geq 2$ locally compact groups. Under an irreducibility assumption, the splitting theorem (Theorem~\ref{thm_splitting}) implies a superrigidity result for uniformly equicontinuous affine $\Gamma$\hy actions on superreflexive spaces $B$. As before such an action can be viewed as an affine isometric $\Gamma$\hy action on
a ucus Banach space $B$. It therefore suffices to apply the splitting theorem to the \emph{induced} $G$\hy action on an \emph{induced space} $L^p(G/\Gamma, B)$ (compare~\cite{Shalom} for the Hilbertian case).

The goal of this section is to address the various (mostly technical) issues that arise when carrying out this programme. We begin by preparing for a statement (Theorem~\ref{thm_super_prod} below) that will then imply a more general form of Theorem~\ref{C:super}.

\subsection{}\label{sec_ind_gen}
Let $G$ be a locally compact group and $\Gamma<G$ a lattice. The induction procedure will work smoothly if $\Gamma$ is uniform (i.e. cocompact); in order to treat some non-uniform cases, one introduces the following.

\begin{defn}[{\cite[III.1.8]{Margulis}}]
The lattice $\Gamma$ is \emph{weakly cocompact} if the $G$\hy representation $L^2_0(G/\Gamma)$, i.e. the canonical complement of the trivial representation in $L^2(G/\Gamma)$, does not almost have non-zero invariant vectors.
\end{defn}

One verifies that any cocompact lattice is weakly cocompact. If $G$ has property~$(T)$, then all its lattices are weakly cocompact. This also holds if $G$ is any connected semisimple Lie group (\cite{Bekka}, compare also~\cite[III.1.12]{Margulis}). By Remark~\ref{rem:T_Lp:specific}, this definition does not depend on considering $L^2(G/\Gamma)$ rather than $L^p(G/\Gamma)$ for some other $1\leq p<\infty$.

\begin{defn}[{See~\cite[1.II]{Shalom}}]\label{defn_p_int}
Let $p>0$. The lattice $\Gamma$ is \emph{$p$\hy integrable} if either~(i) it is uniform; or (ii)~it is finitely generated and for some (or equivalently any) finite generating set $S\se \Gamma$, there is a Borel fundamental domain $\mD\se G$ (with null boundary) such that
$$\int_\mD \|\chi(g^{-1} h)\|_S^p\, d h \ <\infty\kern1cm\forall\,g\in G,$$
where $\|\cdot\|_S$ is the word-length associated to $S$ and $\chi:G\to \Gamma$ is defined by $\chi^{-1}(e)=\mD$, $\chi(g\gamma^{-1})=\gamma\chi(g)$.
\end{defn}

This formulation is a bit awkward so as to include all uniform lattices since~(ii) would otherwise fail when $G$ is not compactly generated. Condition~(ii) holds (with any $p\geq 1$) for all irreducible lattices in higher rank semisimple Lie/algebraic groups, see~\cite[\S2]{Shalom}; it holds likewise for R\'emy's Kac--Moody lattices~\cite{Remy}.

Finally, given a product structure $G=G_1\times \cdots \times G_n$, we say that a lattice $\Gamma<G$ is \emph{irreducible} if its projection to each $G_i$ is dense.

\begin{thm}\label{thm_super_prod}
Let $\Gamma$ be an irreducible lattice in a locally compact $\sigma$\hy compact group $G=G_1\times \cdots \times G_n$. Assume that $\Gamma$ is weakly cocompact and $p$\hy integrable for some $p>1$. Let $B$ be a ucus Banach space with a $\Gamma$\hy action by affine isometries.

If the associated linear $\Gamma$\hy representation does not almost have invariant vectors, then there is a $\Gamma$\hy closed complemented affine subspace of $B$ on which the $\Gamma$\hy action is a sum of actions extending continuously to $G$ and factoring through $G\to G_i$. (Compare Remark~\ref{rem:C:super}.)
\end{thm}

Theorem~\ref{thm_super_prod} indeed implies Theorem~\ref{C:super} in the wider generality of weakly cocompact $p$-integrable lattices, since Proposition~\ref{P:superreflex-affine} allows us to assume that the topological vector space of Theorem~\ref{C:super} is in fact a ucus Banach space with a $\Gamma$\hy action by affine isometries.

\medskip

A (simpler) application of the same techniques implies the following result:

\begin{thm}\label{thm:characters}
Let $\Gamma$ be an irreducible lattice in a locally compact $\sigma$\hy compact group $G=G_1\times \cdots \times G_n$. Assume that $\Gamma$ is weakly cocompact and $p$\hy integrable for some $p>1$.

Then any homomorphism $\Gamma\to\BR$ extends continuously to $G$.
\end{thm}

This result was established by Shalom in the case of cocompact lattices~\cite[0.8]{Shalom} (actually, his proof holds in the setting of square-integrable lattices). It is therefore unsurprising that our reults imply the generalisation stated in Theorem~\ref{thm:characters} above (see the end of this section).

\subsection{Induction}\label{sec_induction}
Throughout this section, $G$ is a locally compact second countable group and $\Gamma<G$ a lattice. In particular, the Haar measure induces a standard Lebesgue space structure on $G/\Gamma$.
\begin{rem}
Even though Theorem~\ref{thm_super_prod} and Theorem~\ref{C:super} was stated in the more general setting of $\sigma$\hy compact groups, it is indeed enough to treat the second countable case: one can reduce to the latter by a structural result of Kakutani--Kodaira~\cite{Kakutani-Kodaira} (the details of the straightforward reduction are expounded at length in~\cite{Monod}).
\end{rem}

Let $B$ be any Banach space and $1<p<\infty$. We consider the Banach space $E=L^p(G/\Gamma, B)$ as in Section~\ref{SS:Lp}.

\begin{lem}\label{lem_ind_uc}
If $B$ is uniformly convex or ucus, then so is $E$.
\end{lem}

\begin{proof}
This follows from a result of Figiel and Pisier, see Theorem~1.e.9 point~(i) in~\cite{Lindenstrauss-Tzafriri}, Volume~II.
\end{proof}

Suppose now that $B$ is endowed with a linear isometric $\Gamma$\hy representation $\ro$. Then $E$ can be canonically isometrically identified
\begin{equation}\label{eq_equiv}
E\ \cong \ L^{[p]}(G,B)^{\ro(\Gamma)}
\end{equation}
with the space of those Bochner-measurable $\Gamma$\hy equivariant function classes $f:G\to B$ such that $\|f\|_B: G/\Gamma\to\BR$ is $p$\hy integrable (the latter condition is symbolized by the notation $L^{[p]}$). Here, we choose to interpret $\Gamma$\hy equivariance as $f(g\gamma)= \ro(\gamma)^{-1}f(g)$. The isomorphism~\eqref{eq_equiv} can be e.g. realized by restricting equivariant maps to any Borel fundamental domain $\mD\se G$ for $\Gamma$ since $\mD\cong G/\Gamma$ as Lebesgue spaces. This identification allows us to endow $E$ with a continuous linear isometric $G$\hy representation by left multiplication. This $G$\hy representation is called the \emph{induced representation}. If we choose a fundamental domain $\mD\se G$ and consider the corresponding map $\chi$ as in Definition~\ref{defn_p_int}, then this $G$\hy representation reads as follows for $f\in E=L^p(G/\Gamma, B)$:
\begin{equation}\label{eq_action_ind}
(hf)(g\Gamma)\ =\ \ro(\chi(g)^{-1}\chi(h^{-1}g)) f(h^{-1}g\chi(h^{-1}g)\Gamma)
\end{equation}
(a good indication that the model~\eqref{eq_equiv} is more natural!).

\begin{lem}\label{lem_ind_aiv}
Assume $\Gamma$ weakly cocompact in $G$. If the linear $\Gamma$\hy representation does not almost have invariant vectors, then the induced linear $G$\hy representation does not either.
\end{lem}

\begin{proof}
The proof given by Margulis in the unitary case~\cite[III.1.11]{Margulis} holds without changes (recalling that we can apply weak cocompactness in the $L^p$ setting by Remark~\ref{rem:T_Lp:specific}).
\end{proof}

Suppose now that $B$ is endowed with an isometric $\Gamma$\hy action~-- not necessarily linear anymore. We want to endow $E$ with a continuous affine isometric $G$\hy action by identifying $E$ with a space of $\Gamma$\hy equivariant function classes $G\to B$ as before, except that equivariance is now understood with respect to the affine $\Gamma$\hy action. Formally, there is nothing to change to the special case of linear action considered above; the action is defined by left $G$\hy translation of equivariant maps, so that via the natural identification we get for $f\in E=L^p(G/\Gamma, B)$ the action
\begin{equation}\label{eq_action_ind:af}
(hf)(g\Gamma)\ =\ \chi(g)^{-1}\chi(h^{-1}g) f(h^{-1}g\chi(h^{-1}g)\Gamma)
\end{equation}
in complete analogy with~\eqref{eq_action_ind}. However, the $L^p$ integrability property might be lost. The condition~(ii) of Definition~\ref{defn_p_int} is a straightforward sufficient condition to retain integrability; cocompactness of $\Gamma$ is also enough, because it ensures that one can choose $\mD$ in such a way that for any compact $C\se G$ the set $\{\eta\in\Gamma : \mD \eta\cap C\neq \varnothing \}$ is finite~\cite[VII \S2 Ex.~12]{Bourbaki}. Compare~\cite[\S2]{Shalom} (and~\cite[App.~B]{Monod}).

In conclusion, we may always consider the continuous \emph{induced (affine) isometric $G$\hy action} on $E$ when $\Gamma$ is $p$\hy integrable.

\smallskip

By construction, the linear part of the induced affine action coincides with the induced linear $G$\hy representation on $E$ considered earlier. If we denote by $b:\Gamma\to B$ the cocycle of the original affine $\Gamma$\hy action, then comparing~\eqref{eq_action_ind} with~\eqref{eq_action_ind:af} shows that the cocycle $\ti b:G\to E$ of the induced affine action is given by
\begin{equation}\label{eq_action_ind:cocycle}
\ti b(h)(g\Gamma)\ =\ b\Big(\chi(g)^{-1}\chi(h^{-1}g)\Big).
\end{equation}
Moreover, the correspondence $b\mapsto \ti b$ induces a (topological) isomorphism $H^1(\Gamma, B)\to H^1(G,E)$.

\smallskip

At this point, we record the following.
\begin{prop}\label{P:lattices}
Keep the notation of this section.
\begin{enumerate}
\item
    If $\Gamma$ has property~$(F_{B})$ then so does $G$.
\item
    If $G$ has property~$(F_E)$ and $\Gamma$ is $p$\hy integrable, then
    $\Gamma$ has property~$(F_{B})$.
\end{enumerate}
\end{prop}
\begin{proof}
For~(1), consider any continuous affine isometric $G$\hy action on $B$; then there is a $\Gamma$\hy fixed point $b\in B$. The corresponding orbit map $G\to B$ descends to a continuous map $G/\Gamma\to B$. The image of the normalized invariant measure on $G/\Gamma$ in $B$ being preserved by $G$, it follows from Lemma~\ref{L:bounded} that there is a $G$\hy fixed point.

For~(2), consider an affine isometric $\Gamma$\hy action on $B$ and endow $E$ with the induced affine action as in the discussion above. Then there is a $G$\hy fixed point $f\in E$. It follows from the description of $E$ as space of equivariant maps that $f$ is essentially constant and that its essential value is a $\Gamma$\hy fixed point of $B$.
\end{proof}

\subsection{Superrigidity}

In order to prove Theorem~\ref{thm_super_prod}, we now analyse the interplay between the induction constructions and the setting of irreducible lattices $\Gamma<G=G_1\times \cdots \times G_n$ as in the beginning of this Section~\ref{S:superrigidity}. We will roughly imitate the arguments given by Shalom in~\cite{Shalom} when he deduces Corollary~4.2 \emph{ibid}.

Keep all the above notations and write $G'_i = \prod_{j\neq i} G_j$. First we observe that the irreducibility of $\Gamma$ implies that for each $i$ it is a well-posed definition to consider the maximal (possibly zero) linear subspace $B_i\se B$ on which the linear $\Gamma$\hy representation $\ro$ extends to a continuous $G$\hy representation $\ro_i:G\to G_i\to \Orth(B_i)$ factoring through $G_i$; moreover $B_i$ is automatically closed by maximality.

The induced space $E$ is ucus by Lemma~\ref{lem_ind_uc}. The isometric (affine) $G$\hy action on $E$ has no fixed point by the very same argument given to prove Proposition~\ref{P:lattices}(2). On the other hand, the linear part does not have almost invariant vectors by Lemma~\ref{lem_ind_aiv}. Thus Theorem~\ref{thm_splitting} applies: There is a $G$\hy invariant closed complemented affine subspace $\underline{E}\se E$ and an affine isometric $G$\hy equivariant isomorphism $\underline{E}\cong E_1 \oplus \cdots \oplus E_n$, where each $E_i$ is a ucus space with an affine isometric $G$\hy action factoring through $G\to G_i$. In view of Remark~\ref{rem:product:subspace}, there is no loss of generality in assuming that $E_i$ is the space of $G'_i$\hy fixed under the induced linear representation. One verifies readily the following:

\begin{lem}\label{lem:induced:fixed}
The map $B_i\to E \cong L^{[p]}(G,B)^{\ro(\Gamma)}$ that to $v\in B_i$ associates the function $G\to B$ defined by $g\mapsto \ro_i(g^{-1}) v$ yields an isometric isomorphism of (linear) $G$\hy spaces $B_i\cong E_i$.\hfill\qedsymbol
\end{lem}

Indeed, since the image of $\Gamma$ in $G_i$ is dense, the Fubini--Lebesgue theorem implies that any map $f:G\to B$ in $E$ that is $G'_i$\hy invariant in the linear representation on $E$ is an orbit map as in the lemma.

\smallskip

At this point we observe that if the subspaces $B_i$ had trivial intersection, we would indeed have found a subspace $\bigoplus B_i \cong \bigoplus E_i$ of $B$ on which the affine $\Gamma$\hy action extends continuously to $G$ as requested. In general, we have a $\Gamma$\hy equivariant affine map
$$\bigoplus E_i \longrightarrow \sum B_i \se B$$
induced by the maps of Lemma~\ref{lem:induced:fixed}. Alternatively, we can think of this map as follows: The cocycle induced as in~\eqref{eq_action_ind:cocycle} decomposes as a sum of cocycles $\ti b=\bigoplus \ti b_i: G\to E$, $\ti b_i:G\to G_i\to E_i$, and in turn by Lemma~\ref{lem:induced:fixed} each $\ti b_i$ is the cocycle induced under the correspondence~\eqref{eq_action_ind:cocycle} from a cocycle $b_i:\Gamma\to B_i$; the affine $\Gamma$\hy action on $\sum B_i$ is determined by the cocycle $\sum b_i$. This completes the proof of Theorem~\ref{thm_super_prod}.\hfill\qedsymbol

\begin{rem}\label{rem:compact:sum}
As mentionned in Remark~\ref{rem:C:super}, the obstruction to extending the affine $\Gamma$\hy action on some subspace of $B$ is confined within a compact group. Indeed, the only reason we might end up with a \emph{sum} of action extending to $G$ through various $G_i$ rather than with a direct sum (which then extends globally to $G$) is the possibility that $B_i\cap B_j\neq 0$ for some $i\neq j$. But then the linear representation of $\Gamma$ on $B_i\cap B_j$ extends continuously to $G$ in two different ways, \emph{both} through $G_i$ and through $G_j$. This may indeed happen but forces the image of $\Gamma$ in $\Orth(B_i\cap B_j)$ to be compact, see examples and discussion in~\cite{Monod}.

Let us only mention the most basic example: $\Gamma< G=G_1\times G_2$ with $G_i = {\bf Z}\rtimes\{\pm1\}$ and $\Gamma={\bf Z}^2\rtimes \{\pm1\}$. Then $\Gamma$ acts affinely isometrically without fixed point on $B=\BR$ (by $(n,m;\e).x = \e x+n+m$) and the associated linear representation does not almost have invariant vectors. However, it is easy to check that this action does not extend to $G$. Instead, it is a sum of actions extending to $G_i$ with sum map $\BR\oplus \BR\to B=\BR$. Here $B_1=B_2=B$.
\end{rem}

\begin{proof}[Proof of Theorem~\ref{thm:characters}]
Recall that the space of homomorphisms $\Gamma\to \BR$ is precisely the space of affine isometric $\Gamma$-actions on $\BR$ with the trivial representation as linear part. By Remark~\ref{rem:T_Lp:specific}, the $G$-representation on $L^p_0(G/\Gamma)$ does not almost have invariant vectors. Therefore, using $p$-integrable induction, one deduces Theorem~\ref{thm:characters} from Theorem~\ref{thm_splitting} very exactly as Shalom deduced Theorem~0.8 in~\cite{Shalom} from Theorem~3.1 in~\cite{Shalom}.
\end{proof}

\section{Appendix: Howe--Moore Theorem on Banach Spaces} \label{Appendix}
In this appendix we sketch the proof of a version of the well known
Howe--Moore theorem on vanishing of matrix coefficients for unitary representations,
extended to the framework of ucus Banach spaces.
This generalization is due to Yehuda Shalom (unpublished) and we state
it here with a sketch of the proof for reader's convenience.
\begin{thm}\label{T:HM}
Let $I$ be a finite set, $k_{i}$, $i\in I$ be local fields,
$\mathbf{G}_{i}$ connected semisimple simply-connected $k_{i}$\hy groups, $G_{i}=\mathbf{G}_{i}(k_{i})$
the locally compact group of $k_{i}$-points, and $G=\prod_{i\in I}
G_{i}$.

Let $B$ be a ucus Banach space and $\ro:G\to \Orth(B)$ a
continuous isometric linear representation,
such that $B^{\ro(G_{i})}=\{0\}$ for each $i\in I$.
Then all matrix coefficients $c_{x,\lambda}(g)=\ip{\ro(g)x}{\lambda}$, $x\in B$, $\lambda\in B^*$,
vanish at infinity, i.e. $c_{x,\lambda}\in C_{0}(G)$.
\end{thm}

Notice that we can (and will) assume that the $\mathbf{G}_{i}$ have no $k_i$-anisotropic factors, since the group of $k_i$-points of such factors are compact.

\begin{proof}[Proof of Theorem~\ref{T:HM}]
In a way of contradiction, assume that for some $g_{n}\to \infty $ in $G$,
$v\in \sphere(B)$, $\lambda\in \sphere(B^{*})$ one has
\[
    \inf |\ip{\ro(g_n)x}{\lambda}|=\gep>0.
\]
We shall prove that at least one simple factor $G_i$ of $G$ has a
non-trivial $\ro(G_i)$ invariant vector.

Let $G=KAK$ be a Cartan decomposition of $G$ (here $K=\prod K_{i}$
and $A=\prod A_{i}$ where $G_{i}=K_{i}A_{i}K_{i}$ is the Cartan
decomposition for $G_i$). We first show that without loss of
generality one may assume $g_{n}\in A$.
\begin{lem}[KAK Reduction] \label{L:KAKreduction}
There exists a sequence $a_n\to\infty$ in the Cartan subgroup $A\se G$ and
non-zero vectors  $y,z\in B$ so that
\[
    \ro(a_n)y\overto{w} z\neq0.
\]
where $\overto{w}$ denotes the weak convergence.
\end{lem}
\begin{proof}
Write $g_n=k_n a_n k^\prime_n$ where $k_n,k^\prime_n\in K$ and $a_n\in A$.
Then $a_{n}\to\infty$ because $g_{n}\to\infty$.
Upon passing to a subsequence, $k^\prime_n\to k^\prime\in K$ and $k_n\to k\in K$.
Denote
\[
    y_n=\ro(k^\prime_n)x,\quad  y=\ro(k^\prime)x,\quad  \mu_n=\ro^*(k_n^{-1})\lambda,
    \quad  \mu=\ro^*(k^{-1})\lambda
\]
where $\ro^*$ is the dual (contragradient) $G$\hy representation on $B^{*}$.
Using the weak-compactness of the unit ball of $B$ we may also assume that
\[
    \ro(a_{n})y\overto{w}z.
\]
We shall show that $\ip{z}{\mu}=\lim \ip{\ro(g_{n})x}{\lambda}$
which is bounded away from zero, hence implying $z\neq 0$.

Recall that in a uc Banach space $B$ the weak and the strong topologies agree on the unit
sphere $\sphere(B)$: indeed if $y_{n}\overto{w} y$ are unit vectors, then
\[
1-\delta(\|y_{n}-y\|) \geq  \|y_{n}+y\|/2 \geq \ip{(y_{n}+y)/2}{y^{*}} \to 1.
\]
Hence $\delta(\|y_{n}-y\|)\to 0$ and $\|y_{n}-y\|\to0$.
For the same reason we also have $\|\mu_n-\mu\|\to 0$ in $\sphere(B^{*})$.
For an arbitrary $\xi\in B^*$
\[
    |\ip{\ro(a_n)y_n}{\xi}-\ip{\ro(a_n)y}{\xi}|\le\|y_n-y\|\cdot\|\xi\|\to0.
\]
Hence $\ro(a_{n})y_{n}\overto{w} z$.
In general, if $z_n\overto{w}z$ in $B$ and $\mu_n\overto{w}\mu$ in $B^*$ then $\ip{z_n}{\mu_n}\to\ip{z}{\mu}$ because weakly convergent sequences are bounded
in norm and
\begin{eqnarray*}
    |\ip{z_n}{\mu_n}- \ip{z}{\mu}|&\le&|\ip{z_n}{\mu_n-\mu}|+|\ip{z_n-z}{\mu}|\\
    &\le& (\sup\|z_n\|)\cdot \|\mu_n-\mu\|^*+|\ip{z_n-z}{\mu}|\to0.
\end{eqnarray*}
Therefore
\[
    \ip{\ro(g_n)x}{\lambda}=\ip{\ro(a_nk^\prime_n)x}{\ro^*(k_n^{-1})\lambda}        =\ip{\ro(a_n)y_n}{\mu_n}\to\ip{z}{\mu}
\]
implying $|\ip{z}{\mu}|\ge \gep$, which in particular means that $z\neq0$.
\end{proof}
\begin{lem}[Generalized Mautner Lemma] \label{L:Mautner}
Suppose that $\{a_n\}$ and $h$ in $G$ satisfy $a_{n}^{-1}h a_{n}\to 1_{G}$ in $G$.
If $y,z\in B$ are such that $\ro(a_n)y\overto{w}z$ then $\ro(h)z=z$.
In particular, if $\ro(a_{n})z=z$ then $\ro(h)z=z$.
\end{lem}
\begin{proof}
(Strong) continuity of $\ro$ gives
\[
    \|\ro(ha_n)y-\ro(a_n)y\|=\left\|\ro(a_n^{-1}ha_n)y-y\right\|\to0
\]
At the same time $\ro(a_{n})y\overto{w} z$ and $\ro(ha_{n})y\overto{w}\ro(h) z$.
Hence $\ro(h)z=z$.
\end{proof}

We can now prove Theorem~\ref{T:HM} in the case of $G=\SL_{2}(k)$
where $k$ is a local field. Assuming that
$\ro:\SL_{2}(k)\to\Orth(B)$ has some matrix coefficient not
vanishing at infinity, we get by Lemma~\ref{L:KAKreduction} a
sequence $a_n\to\infty$ in  $A$,
and \emph{non zero}
vectors $y,z\in B$ with $\ro(a_n)y\overto{w} z$.

Let $H$ be the unipotent (horocyclic) subgroup defined by
$H=\{h\in G:a_n^{-1}ha_n\to e\}$. It is normalized by $a_n$, and
by Lemma~\ref{L:Mautner} $z$ is a (non-trivial) $\ro(H)$\hy
invariant vector. We may assume that $\|z\|=1$. The matrix
coefficient $f(g)=\ip{\ro(g)z}{z^{*}}$ is a continuous function on
$G$, which is bi-$H$\hy invariant:
\begin{align}
\label{former_a}
f(g h)&=\ip{\ro(g)\ro(h)z}{z^{*}}=\ip{\ro(g)z}{z^{*}}=f(g)\\
\label{former_b}
f(h g)&=\ip{\ro(g)z}{\ro^*(h^{-1})z^{*}}=\ip{\ro(g)z}{z^{*}}=f(g)
\end{align}
for all $g\in G$ and $h\in H$. The proof can be now completed as
in the original unitary Howe--Moore Theorem. By~\eqref{former_a}, $f$ can be
viewed as a continuous function $f_0$ on the punctured plane
$G/H=k^2-\{(0,0)\}$, and by~\eqref{former_b}, $f_0$ is constant on each
horizontal line $\ell_s=\{(t,s)\,:\, t\in k\}$, $s\neq 0$, where
we identify $H$ with the upper triangular unipotent subgroup by
choosing an appropriate basis for $k^2$. By continuity, $f_0$ is a
constant on $\{(t,0)\,:\,t\neq 0\}$. Since $f_0(0,1)=f(e)=1$ this
constant is $1$.

This implies that $z$ is $\ro (A)$-invariant because
$\ip{\ro(a)z}{z^{*}}=f(a)=f(e)=1$ whilst $z^{*}$ attains its norm
only on $z$.

Thus $z$ is fixed by the upper triangular group $AH\se G$
and $f$ descends to
a continuous function $f_1$ on the projective line $\BP(k^{2})=G/AH$.
The $H$\hy action on $\BP(k^{2})$ has a dense orbit.
Thus $f_{1}$ is constant $1$, and so is $f$:
\[
    \ip{\ro(g)z}{z^{*}}=f(g)=f(e)=1\qquad (g\in G)
\]
Thus the unit vector $z$ is $\ro(G)$\hy invariant, completing the proof in
the case of $G=\SL_{2}(k)$.

The proof of the unitary Howe--Moore theorem for semisimple Lie group $G=\prod G_i$
(c.f. Zimmer~\cite{Zimmer84}, Margulis~\cite{Margulis}) relies only on the reduction to the
Cartan subgroup (Lemma~\ref{L:KAKreduction}),
the structure of such groups, the case of $\SL_{2}(k)$
and on Mautner Lemma.
Thus the ``unitary'' argument can be applied almost \emph{verbatim} to the present setup of
ucus Banach spaces.
\end{proof}


\end{document}